# ON THE CONVERGENCE OF THE ZETA FUNCTION FOR CERTAIN PREHOMOGENEOUS VECTOR SPACES

Akihiko Yukie[1]

Oklahoma State University

**Introduction**

Let $(G, V)$ be an irreducible prehomogeneous vector space defined over a number field $k$, $P \in k[V]$ a relative invariant polynomial, and $\chi$ a rational character of $G$ such that $P(gx) = \chi(g)P(x)$. Let $V_k^{\text{ss}} = \{x \in V_k \mid P(x) \neq 0\}$. For $x \in V_k^{\text{ss}}$, let $G_x$ be the stabilizer of $x$, and $G_x^0$ the connected component of 1 of $G_x$. We define $L_0$ to be the set of $x \in V_k^{\text{ss}}$ such that $G_x^0$ does not have a non-trivial rational character. Then we define the zeta function for $(G, V)$ by the following integral

$$Z(\Phi, s) = \int_{G_{\mathbb{A}}/G_k} |\chi(g)|^s \sum_{x \in L_0} \Phi(gx) dg,$$

where $\Phi$ is a Schwartz–Bruhat function, $s$ is a complex variable, and $dg$ is an invariant measure.

Shintani showed the convergence of $Z(\Phi, s)$ for $\text{Re}(s) \gg 0$ for the spaces $\text{Sym}^2 k^n$ and $\text{Sym}^3 k^2$ (see [4], [5]). F. Sato showed the convergence of $Z(\Phi, s)$ when $G_x \cap \text{Ker}(\chi)$ is connected semi-simple (which implies that $L_0 = V_k^{\text{ss}}$) (see [1]). Note that his assumptions in [1] were later proved by other people. Also he considered prehomogeneous vector spaces over $\mathbb{Q}$, but if $(G, V)$ is a prehomogeneous vector space over $k$, we can consider $(G, V)$ as a prehomogeneous vector space over $\mathbb{Q}$. Then the zeta function of $(G, V)$ over $k$ and the zeta function of $(G, V)$ over $\mathbb{Q}$ are the same. So his result implies the convergence of the zeta function for prehomogeneous vector spaces as above over an arbitrary number field $k$. In [8], we showed the convergence of $Z(\Phi, s)$ when $\dim G = \dim V$ (in this case $L_0 = V_k^{\text{ss}}$ also). These cover 23 types of irreducible reduced prehomogeneous vector spaces. Ying recently showed the convergence of $Z(\Phi, s)$ for a few cases when $L_0 \neq V_k^{\text{ss}}$. In this paper, we prove the convergence of $Z(\Phi, s)$ for prehomogeneous vector spaces of the form $(G/\widetilde{T}, V)$, where $G, V$ are as follows:
  (1) $G = \text{GL}(2) \times \text{GL}(2) \times \text{GL}(2)$, $V = k^2 \otimes k^2 \otimes k^2$,
  (2) $G = \text{GL}(3) \times \text{GL}(3) \times \text{GL}(2)$, $V = k^3 \otimes k^3 \otimes k^2$,
  (3) $G = \text{GL}(4) \times \text{GL}(2)$, $V = \wedge^2 k^4 \otimes k^2$,
  (4) $G = \text{GL}(6) \times \text{GL}(2)$, $V = \wedge^2 k^6 \otimes k^2$,

[1]Partially supported by NSF grant DMS-9401391 and NSA grant MDA-904-93-H-3035

and $\widetilde{T} = \mathrm{Ker}(G \to \mathrm{GL}(V))$ for all the cases. These are the $D_4$, $E_6$, $D_5$, $E_7$ cases in [6].

Note that since $L_0 = V_k^{\mathrm{ss}}$ for the case (4), the result of M. Sato and Shintani (see [3]) on the meromorphic continuation and the functional equation of the local zeta function at an infinite place implies the meromorphic continuation of $Z(\Phi, s)$ and the functional equation of the form
$$Z(\Phi, s) = Z(\widehat{\Phi}, N - s),$$
where $\widehat{\Phi}$ is an appropriate Fourier transform and $N$ is a number which can easily be figured out depending on the normalization. (see §0.3 of [8]). For the cases (1)–(3), the meromorphic continuation of $Z(\Phi, s)$ is unknown.

In [7], Ying considered three types of prehomogeneous vector spaces, one of which is the case where $G = \mathrm{GSpin}(Q) \times \mathrm{GL}(2)$ for a non-degenerate quadratic form $Q$ in $n \geq 4$ variables, and $V$ is the tensor product of the standard representations. When $\mathrm{GSpin}(Q)$ is split, the case $n = 4$ (resp. $n = 6$) is the case (1) (resp. case (3)) of this paper. So cases (1) and (3) of this paper are covered by Ying. However, our method is totally different from Ying's method. For example, his method is based on the consideration of Tamagawa numbers as in F. Sato's paper [1] and does not prove that the incomplete theta series $\sum_{x \in L_0} \Phi(gx)$ satisfies the assumption of Shintani's lemma (see §3.4 of [8]). Our method is to estimate the incomplete theta series on a Siegel set. Therefore, we can show that $\sum_{x \in L_0} \Phi(gx)$ satisfies the assumption of Shintani's lemma.

We handle the cases (1), (2) in §2, and the cases (3), (4) in §3.

## §1 Preliminaries

We basically follow the notations of [8], but we recall the most basic ones. For a finite set $X$, $\#X$ is its cardinality. If $f, g$ are functions on a set $X$ (not necessarily finite), $f \ll g$ means that there exists a constant $C$ such that $f(x) \leq Cg(x)$ for all $x \in X$. We also use the classical notation $x \ll y$ when $y$ is a much larger number than $x$. We hope the meaning of this notation will be clear from the context. The ring of adeles (resp. the group of ideles) over $k$ is denoted by $\mathbb{A}$ (resp. $\mathbb{A}^\times$). For a vector space $V$ over $k$, $V_\mathbb{A}$ is the adelization, and $\mathscr{S}(V_\mathbb{A})$ is the space of Schwartz–Bruhat functions. We define $\mathbb{R}_+ = \{x \in \mathbb{R} \mid x > 0\}$. For $\lambda \in \mathbb{R}_+$, $\underline{\lambda}$ is the idele whose component at any infinite place is $\lambda^{\frac{1}{[k:\mathbb{Q}]}}$ and whose component at any finite place is 1. Let $|x|$ be the adelic absolute value of $x \in \mathbb{A}$. Then $|\underline{\lambda}| = \lambda$. Let $a_n(t_1, \cdots, t_n)$ be the $n$-dimensional diagonal matrix whose $(i, i)$-entry is $t_i$ for all $i$. We define $\mathrm{GL}(n)_\mathbb{A}^0 = \{g \in \mathrm{GL}(n)_\mathbb{A} \mid |\det g| = 1\}$.

For all the four cases in this paper, $G$ is of the form $G = \mathrm{GL}(n_1) \times \cdots \mathrm{GL}(n_f)$. ($f$ is either 2 or 3). Let $G_i = \mathrm{GL}(n_i)$ for all $i$. Let $T_i \subset G_i$ be the set of diagonal matrices, and $T = T_1 \times \cdots T_f$. Let $\epsilon > 0$ be a sufficiently small constant. We define
$$G_{i\mathbb{A}}^0 = \mathrm{GL}(n_i)_\mathbb{A}^0,$$
$$T_{i+}^0 = \{a_{n_i}(\underline{\lambda}_{i1}, \cdots, \underline{\lambda}_{in_i}) \mid \lambda_{i1}, \cdots, \lambda_{in_i} \in \mathbb{R}_+, \ \lambda_{i1} \cdots \lambda_{in_i} = 1\},$$
$$T_{i\epsilon}^0 = \left\{ a_{n_i}(\underline{\lambda}_{i1}, \cdots, \underline{\lambda}_{in_i}) \left| \begin{array}{l} \lambda_{i1}, \cdots, \lambda_{in_i} \in \mathbb{R}_+, \ \lambda_{i1} \cdots \lambda_{in_i} = 1 \\ \lambda_{i1}\lambda_{i2}^{-1}, \cdots, \lambda_{in_i-1}\lambda_{in_i}^{-1} \geq \epsilon \end{array} \right. \right\},$$
$$\mathfrak{t}_i^* = \{y_i = (y_{i1}, \cdots, y_{in_i}) \in \mathbb{R}^{n_i} \mid y_{i1} + \cdots + y_{in_i} = 0\}.$$

For $c_i = (c_{i1}, \cdots, c_{in_i-1}) \in \mathbb{R}^{n_i-1}$, we define
$$w_i(c_i) = c_{i1}(1, -1, 0, \cdots, 0) + c_{i2}(0, 1, -1, 0, \cdots, 0) + \cdots + c_{in_i-1}(0, \cdots, 0, 1, -1).$$

Let $\mathfrak{t}^*_{i,\mathrm{pc}}$ be the cone generated by positive weights, i.e.
$$\mathfrak{t}^*_{i,\mathrm{pc}} = \{w_i(c_i) \in \mathfrak{t}^*_i \mid c_{i1}, \cdots, c_{in_i-1} \geq 0\}.$$

Apparently, the set of interior points of $\mathfrak{t}^*_{i,\mathrm{pc}}$ consists of points of the form $w_i(c_i)$ where $c_{i1}, \cdots, c_{in_i-1} > 0$. For $c = (c_1, \cdots, c_f)$, we define
$$w(c) = (w_1(c_1), \cdots, w_f(c_f)).$$

Let
$$G^0_{\mathbb{A}} = G^0_{1\mathbb{A}} \times \cdots \times G^0_{f\mathbb{A}},$$
$$T^0_+ = T^0_{1+} \times \cdots \times T^0_{f+},$$
$$T^0_\epsilon = T^0_{1\epsilon} \times \cdots \times T^0_{f\epsilon},$$
$$\mathfrak{t}^* = \mathfrak{t}^*_1 \times \cdots \times \mathfrak{t}^*_f,$$
$$\mathfrak{t}^*_{\mathrm{pc}} = \mathfrak{t}^*_{1,\mathrm{pc}} \times \cdots \times \mathfrak{t}^*_{f,\mathrm{pc}}.$$

For $t \in T^0_+$ and $y = (y_1, \cdots, y_f) \in \mathfrak{t}^*$, we can define $t^y \in \mathbb{R}_+$ in the usual manner. Let $\rho \in \mathfrak{t}^*$ be half the sum of positive weights. This means that $t^\rho = \prod_i \prod_{j<k} (\lambda_{ij}\lambda_{ik}^{-1})$ for
$$t = (a_{n_1}(\underline{\lambda}_{11}, \cdots, \underline{\lambda}_{1n_1}), \cdots, a_{n_f}(\underline{\lambda}_{11}, \cdots, \underline{\lambda}_{fn_f})).$$

The Weyl group $W$ of $G$ is the product of the Weyl groups of $\mathrm{GL}(n_1), \cdots, \mathrm{GL}(n_f)$ and we identify the Weyl group of $\mathrm{GL}(n_i)$ as the set of permutation matrices for all $i$. The group $W$ acts on $T^0_+$ from the left by $t \to gtg^{-1}$ for $g \in W$, $t \in T^0_+$. We define the left action of $W$ on $\mathfrak{t}^*$ by $t^{gy} = (g^{-1}tg)^y$ for $g \in W$, $y \in \mathfrak{t}^*$, $t \in T^0_+$.

For the cases in this paper, up to a constant, $Z(\Phi, s)$ coincides with the following integral
$$\int_{\mathbb{R}_+ \times G^0_{\mathbb{A}}/G_k} \underline{\lambda}^s \sum_{x \in L_0} \Phi(\underline{\lambda}g^0 x) d^\times \underline{\lambda} dg^0,$$
where the action of $\underline{\lambda}$ is the usual multiplication by $\underline{\lambda}$, $d^\times \lambda = \lambda^{-1} d\lambda$, and $dg^0$ is an invariant measure on $G^0_{\mathbb{A}}$. We define
$$Z_+(\Phi, s) = \int_{[1,\infty) \times G^0_{\mathbb{A}}/G_k} \underline{\lambda}^s \sum_{x \in L_0} \Phi(\underline{\lambda}g^0 x) d^\times \underline{\lambda} dg^0.$$

It is well known that there exists a compact set $\widehat{\Omega} \subset G^0_{\mathbb{A}}$ such that $\widehat{\Omega} T^0_\epsilon$ surjects to $G^0_{\mathbb{A}}/G_k$. Therefore, by Proposition (1.2.3) [8], there exists $0 \leq \Psi \in \mathscr{S}(V_{\mathbb{A}})$ such that $Z(\Phi, s), Z_+(\Phi, s)$ are bounded by constant multiples of the following integrals

(1.1)
$$\int_{\mathbb{R}_+ \times T^0_\epsilon} \underline{\lambda}^{\mathrm{Re}(s)} \sum_{x \in L_0} \Psi(\underline{\lambda} tx) t^{-2\rho} d^\times \underline{\lambda} d^\times t,$$
$$\int_{[1,\infty) \times T^0_\epsilon} \underline{\lambda}^{\mathrm{Re}(s)} \sum_{x \in L_0} \Psi(\underline{\lambda} tx) t^{-2\rho} d^\times \underline{\lambda} d^\times t$$

respectively, where $d^\times t$ is an invariant measure on $T_+^0$.

In the following sections, we choose a coordinate system $x = (x_1, \cdots, x_N)$ of $V$ for each case so that there exists $\gamma_i \in \mathfrak{t}^*$ for $i = 1, \cdots, N$ and $tx = (\underline{t^{\gamma_i} x_i})$ for $t \in T_+^0$, $x \in V_\mathbb{A}$. The element $\gamma_i$ is called the weight of the coordinate $x_i$. For $x = (x_1, \cdots, x_N) \in V_k$, we define $I_x = \{1 \leq i \leq N \mid x_i \neq 0\}$. Let $\mathrm{Conv}_x$ be the convex hull of the set $\{\gamma_i \mid i \in I_x\}$.

**Definition (1.2)** *A point $x \in V_k$ is $k$-stable if for all $g \in G_k$, the convex hull $\mathrm{Conv}_{gx}$ contains a neighborhood of the origin of $\mathfrak{t}^*$.*

We showed in Proposition (3.1.4) [8] that if $L_0$ coincides with the set of $k$-stable points, $Z(\Phi, s)$ converges absolutely for $\mathrm{Re}(s) \gg 0$ and $Z_+(\Phi, s)$ is an entire function.

We need the following lemma in §2 to show that $L_0$ coincides with the set of $k$-stable points for the cases (1), (2).

**Lemma (1.3)** *Suppose that $L \subset V_k^{\mathrm{ss}}$ is a $G_k$-invariant subset such that $\mathrm{Conv}_x$ contains an interior point of $\mathfrak{t}_{\mathrm{pc}}^*$ for any $x \in L$. Then $x$ is $k$-stable for all $x \in L$.*

*Proof.* Suppose $x \in L$. Let $g \in W$, $t \in T_+^0$. We define

$$e_i = (\overbrace{0, \cdots, 0, 1}^{i}, 0, \cdots, 0) \in V_k$$

for $i = 1, \cdots, N$. Then $te_i = \underline{t^{\gamma_i}} e_i$. So

$$tge_i = gg^{-1}tge_i = g\underline{(g^{-1}tg)^{\gamma_i}} e_i = \underline{t^{g\gamma_i}} ge_i.$$

Therefore, $\mathrm{Conv} I_{gx} = g\mathrm{Conv}_x$. Since $L_0$ is $G_k$-invariant, $gx \in L_0$. This implies that $g\mathrm{Conv}_x$ contains an interior point of $\mathfrak{t}_{\mathrm{pc}}^*$. So $\mathrm{Conv}_x$ contains an interior point of $g^{-1}\mathfrak{t}_{\mathrm{pc}}^*$. Note that this statement is true for all $g \in W$.

Suppose that $\mathrm{Conv}_x$ does not contain a neighborhood of the origin of $\mathfrak{t}^*$. Since $\mathrm{Conv}_x$ is a finite convex polytope, this implies that $\mathrm{Conv}_x$ is contained in a half space containing the origin, say $\{y \in \mathfrak{t}^* \mid l(y) \leq 0\}$ where $l(y)$ is a non-zero linear form on $\mathfrak{t}^*$. There exists an element $g \in W$ such that $l(g^{-1}y)$ is of the form

$$l(g^{-1}y) = l_1(y_1) + \cdots + l_f(y_f),$$
$$l_i(y_i) = a_{i1}y_{i1} + \cdots + a_{in_i}y_{in_i}$$

for $y = (y_1, \cdots, y_f) \in \mathfrak{t}^*$ where $a_{i1} \geq \cdots \geq a_{in_i}$ are constants for $i = 1, \cdots, f$.

Since $y_{i1} + \cdots + y_{in_i} = 0$ for all $i$, we may assume that $a_{ij} > 0$ for all $i, j$. Also since the linear form $l$ is not identically zero, we may assume that there exist $i_0, j_0$ such that $a_{i_0 j_0} > a_{i_0 j_0+1}$.

We showed that there exists an interior point $w(c) = (w_1(c_1), \cdots, w_f(c_f))$ of $\mathfrak{t}_{\mathrm{pc}}^*$ such that $\mathrm{Conv}_x$ contains the point $g^{-1}w(c)$. Then

$$l(g^{-1}w(c)) = \sum_{i=1}^{f} \sum_{j=1}^{n_i-1} (a_{ij} - a_{ij+1})c_{ij}.$$

By assumption, all the terms are non-negative and at least one term is positive. Therefore, $l(g^{-1}w(c)) > 0$. This is a contradiction. So we can conclude that $\text{Conv}_x$ contains a neighborhood of the origin.

Q.E.D.

§**2** $D_4$, $E_6$ **cases**

We consider the cases (1), (2) in the introduction in this section. We consider these prehomogeneous vector spaces as $\text{M}(2,2) \otimes k^2$ or $\text{M}(3,3) \otimes k^2$, i.e. the space of $2 \times 2$ or $3 \times 3$ matrices whose entries are linear forms in two variables $v = (v_1, v_2)$. We express a general element of $V$ as $M_x(v) = v_1 x_1 + v_2 x_2$ where $x_1 = (x_{1,ij})$, $x_2 = (x_{2,ij})$ are $2 \times 2$ or $3 \times 3$ matrices. We choose $x = (x_1, x_2)$ as the coordinate system of $V$. If $g = (g_1, g_2, g_3)$ is an element of $\text{GL}(2) \times \text{GL}(2) \times \text{GL}(2)$ or $\text{GL}(3) \times \text{GL}(3) \times \text{GL}(2)$, the action of $g$ is defined by

$$gM_x(v) = g_1 M_x(vg_3)^t g_2.$$

We define $F_x(v) = \det M_x(v)$. Then $F_x$ is a binary quadratic or cubic form. It was proved in [2] that $V_k^{\text{ss}}$ is the set of $x$ such that $F_x$ has distinct factors over the closure $\bar{k}$ of $k$. We showed in [6] that $L_0$ is the set of $x$ such that $F_x$ is irreducible.

**Theorem (2.1)** *The set $L_0$ coincides with the set of $k$-stable points. Therefore, $Z(\Phi, s)$ converges absolutely for $\text{Re}(s) \gg 0$ and $Z_+(\Phi, s)$ is an entire function.*

*Proof.* Suppose that $F_x$ is irreducible. Then for any $v \in k^2 \setminus \{(0,0)\}$, $F_x(v) \neq 0$, i.e. $M_x(v)$ is a non-singular matrix. In particular $x_1, x_2$ are non-singular matrices. Let $t = (t_1, t_2, t_3) \in T_+^0$, where

$$\begin{cases} t_1 = a_2(\underline{\lambda}_{11}, \underline{\lambda}_{11}^{-1}), \; t_2 = a_2(\underline{\lambda}_{21}, \underline{\lambda}_{21}^{-1}), \; t_3 = a_2(\underline{\lambda}_{31}, \underline{\lambda}_{31}^{-1}) & \text{case (1)}, \\ t_1 = a_3(\underline{\lambda}_{11}, \underline{\lambda}_{12}, \underline{\lambda}_{13}), \; t_2 = a_3(\underline{\lambda}_{21}, \underline{\lambda}_{22}, \underline{\lambda}_{23}), \; t_3 = a_2(\underline{\lambda}_{31}, \underline{\lambda}_{31}^{-1}) & \text{case (2)}, \end{cases}$$

and $\lambda_{11}\lambda_{12}\lambda_{13} = \lambda_{21}\lambda_{22}\lambda_{23} = 1$.

The set $L_0$ is clearly $G_k$-invariant. So by Lemma (1.3), we only have to show that for any $x \in L_0$, $\text{Conv}_x$ contains an interior point of $\mathfrak{t}_{\text{pc}}^*$. Let $\gamma_{i,jk}$ be the weight of the coordinate $x_{i,jk}$ for all $i, j, k$. The element $\gamma_{i,jk}$ can be expressed in the form $\gamma_{i,jk} = w(d_{i,jk})$ ($d_{i,jk}$ may not be in $\mathfrak{t}_{\text{pc}}^*$).

We first consider the case (1). The following lemma is easy to verify and the proof is left to the reader.

**Lemma (2.2)** (1) $d_{1,11} = (\frac{1}{2}, \frac{1}{2}, \frac{1}{2})$.
(2) $d_{1,12} = (\frac{1}{2}, -\frac{1}{2}, \frac{1}{2})$.
(3) $d_{1,21} = (-\frac{1}{2}, \frac{1}{2}, \frac{1}{2})$.

Note that $d_{2,jk}$ can be obtained by replacing the last $\frac{1}{2}$ in $d_{1,jk}$ by $-\frac{1}{2}$.

Suppose $x \in L_0$. If $x_{1,11} \neq 0$, then $\gamma_{1,11} \in \text{Conv}_x$ and $\gamma_{1,11}$ is an interior point of $\mathfrak{t}_{\text{pc}}^*$ by the above lemma.

Suppose $x_{1,11} = 0$. Then since $x_1$ is non-singular, $x_{1,12}, x_{1,21} \neq 0$. Moreover if $x_{2,11} = 0$, we can choose $v \in k^2 \setminus \{0\}$ so that $M_x(v)$ is singular. This contradicts to the assumption $x \in L_0$. So we may assume that $x_{2,11} \neq 0$. Therefore, $\gamma_{1,12}, \gamma_{1,21}, \gamma_{2,11} \in \text{Conv}_x$. This implies that $\gamma_{1,12} + \gamma_{1,21} + \gamma_{2,11} \in \text{Conv}_x$ also and

$$d_{1,12} + d_{1,21} + d_{2,11} = \left(\frac{1}{2}, \frac{1}{2}, \frac{1}{2}\right).$$

So $\gamma_{1,12} + \gamma_{1,21} + \gamma_{2,11}$ is an interior point of $\text{Conv}_x$. This completes the proof of Theorem (2.1) for the case (1).

Next, we consider the case (2). The following lemma is easy to verify and the proof is left to the reader.

**Lemma (2.3)** (1) $d_{1,11} = ((\frac{2}{3}, \frac{1}{3}), (\frac{2}{3}, \frac{1}{3}), \frac{1}{2})$.
(2) $d_{1,12} = ((\frac{2}{3}, \frac{1}{3}), (-\frac{1}{3}, \frac{1}{3}), \frac{1}{2})$.
(3) $d_{1,13} = ((\frac{2}{3}, \frac{1}{3}), (-\frac{1}{3}, -\frac{2}{3}), \frac{1}{2})$.
(4) $d_{1,21} = ((-\frac{1}{3}, \frac{1}{3}), (\frac{2}{3}, \frac{1}{3}), \frac{1}{2})$.
(5) $d_{1,22} = ((-\frac{1}{3}, \frac{1}{3}), (-\frac{1}{3}, \frac{2}{3}), \frac{1}{2})$.
(6) $d_{1,31} = ((-\frac{1}{3}, -\frac{2}{3}), (\frac{2}{3}, \frac{1}{3}), \frac{1}{2})$.

Note that $d_{2,jk}$ can be obtained by replacing the last $\frac{1}{2}$ in $d_{1,jk}$ by $-\frac{1}{2}$.

Suppose $x \in L_0$. If $x_{1,11} \neq 0$, then $\gamma_{1,11} \in \text{Conv}_x$ and $\gamma_{1,11}$ is an interior point of $\mathfrak{t}_{\text{pc}}^*$.

Suppose $x_{1,11} = 0$, $x_{1,12}, x_{1,21} \neq 0$. Then $\gamma_{1,12}, \gamma_{1,21} \in \text{Conv}_x$. So $\gamma_{1,12} + \gamma_{1,21} \in \text{Conv}_x$ also and

$$d_{1,12} + d_{1,21} = \left(\left(\frac{1}{3}, \frac{2}{3}\right), \left(\frac{1}{3}, \frac{2}{3}\right), 1\right).$$

So $\gamma_{1,12} + \gamma_{1,21}$ is an interior point of $\mathfrak{t}_{\text{pc}}^*$.

Consider the following two cases:
(1) $x_{1,11} = 0$, $x_{1,12} = 0$, and $x_{1,21} \neq 0$,
(2) $x_{1,11} = 0$, $x_{1,21} = 0$, and $x_{1,12} \neq 0$.

Since these cases are similar, we only consider the case (1). Since $x_1$ is a non-singular matrix, $x_{1,13} \neq 0$. If $x_{2,11} = x_{2,12} = 0$, we can choose $v \in k^2 \setminus \{0\}$ so that $v_1 x_{1,13} + v_2 x_{2,13} = 0$. This contradicts to the assumption $x \in L_0$. So we may assume that either $x_{2,11} \neq 0$ or $x_{2,12} \neq 0$. Since

$$d_{2,11} = d_{2,12} + ((0,0), (1,0), 0),$$

we only consider the case $x_{2,12} \neq 0$.

With these assumptions, $\gamma_{1,21}, \gamma_{1,13}, \gamma_{2,12} \in \text{Conv}_x$. Then

$$3\gamma_{1,21} + 2\gamma_{1,13} + 2\gamma_{2,12} \in \text{Conv}_x$$

also and

$$3d_{1,21} + 2d_{1,13} + 2d_{2,12} = \left(\left(\frac{5}{3}, \frac{7}{3}\right), \left(\frac{2}{3}, \frac{1}{3}\right), \frac{3}{2}\right).$$

So $3\gamma_{1,21} + 2\gamma_{1,13} + 2\gamma_{2,12}$ is an interior point of $\mathfrak{t}_{\text{pc}}^*$.

Suppose $x_{1,11} = x_{1,12} = x_{1,21} = 0$. Then since $x_1$ is a non-singular matrix, $x_{1,13}, x_{1,22}, x_{1,31} \neq 0$. Suppose $x_{2,11} \neq 0$. Then

$$\gamma_{1,13} + \gamma_{1,22} + \gamma_{1,31} + \gamma_{2,11} \in \text{Conv}_x$$

and
$$d_{1,13} + d_{1,22} + d_{1,31} + d_{2,11} = \left(\left(\frac{2}{3}, \frac{1}{3}\right), \left(\frac{2}{3}, \frac{1}{3}\right), 1\right).$$

So $\gamma_{1,13} + \gamma_{1,22} + \gamma_{1,31} + \gamma_{2,11}$ is an interior point of $\mathfrak{t}_{pc}^*$.

Suppose $x_{1,11} = x_{1,12} = x_{1,21} = x_{2,11} = 0$. Then if either $x_{2,12} = 0$ or $x_{2,21} = 0$, we can choose $v \in k^2 \setminus \{0\}$ so that $v_1 x_1 + v_2 x_2$ is singular, which is a contradiction. So $x_{2,12}, x_{2,21} \neq 0$. By assumption, $x_{1,13}, x_{1,22}, x_{1,31} \neq 0$ also. Then

$$\gamma_{1,13} + \gamma_{1,22} + \gamma_{1,31} + \gamma_{2,12} + \gamma_{2,21} \in \text{Conv}_x$$

and
$$d_{1,13} + d_{1,22} + d_{1,31} + d_{2,12} + d_{2,21} = \left(\left(\frac{1}{3}, \frac{2}{3}\right), \left(\frac{1}{3}, \frac{2}{3}\right), \frac{1}{2}\right).$$

So $\gamma_{1,13} + \gamma_{1,22} + \gamma_{1,31} + \gamma_{2,12} + \gamma_{2,21}$ is an interior point of $\mathfrak{t}_{pc}^*$. This completes the proof of Theorem (2.1) for the case (2).

Q.E.D.

## §3 $D_5$, $E_7$ cases

We consider the cases (3), (4) in the introduction in this section. We consider these cases as the space of $4 \times 4$ or $6 \times 6$ alternating matrices whose entries are linear forms in two variables $v = (v_1, v_2)$. We express a general element of $V$ as $M_x(v) = v_1 x_1 + v_2 x_2$ where $x_1 = (x_{1,ij})$, $x_2 = (x_{2,ij})$ are $4 \times 4$ or $6 \times 6$ alternating matrices. We choose $x = (x_1, x_2)$ as the coordinate system of $V$ (we only consider $x_{i,jk}$ such that $j > k$). If $g = (g_1, g_2)$ is an element of $\text{GL}(4) \times \text{GL}(2)$ or $\text{GL}(6) \times \text{GL}(2)$, the action of $g$ is defined by

$$gM_x(v) = g_1 M_x(vg_2)^t g_1.$$

Since $M_x(v)$ is an alternating matrix, there exists a binary quadratic or cubic form $F_x(v)$ such that $\det M_x(v) = F_x(v)^2$ ($F_x(v)$ is the Pfaffian of $M_x(v)$). It was proved in [2] that $V_k^{ss}$ is the set of $x$ such that $F_x(v)$ has distinct factors. We showed in [6] that $L_0$ is the set of $x$ such that $F_x$ is irreducible for the case (3) and that $L_0 = V_k^{ss}$ for the case (4).

**Theorem (3.1)** *The integral $Z(\Phi, s)$ converges absolutely and locally uniformly for $\text{Re}(s) \gg 0$ and $Z_+(\Phi, s)$ is an entire function.*

*Proof.* Unlike the cases (1), (2), there are no $k$-stable points, so we have to be a little more subtle for these cases. Let $\Psi$ be as in §1. For $L \subset V_k$, we define

(3.2) $$\Theta_L(\Psi, \underline{\lambda} t) = \sum_{x \in L} \Psi(\underline{\lambda} t x)$$

for $\lambda \in \mathbb{R}_+$, $t \in T_\epsilon^0$.

We estimate $\Theta_{L_0}(\Psi, \underline{\lambda} t)$. Note that if $y \in \mathfrak{t}^*$, the integral $\int_{T_\epsilon^0} t^{y-2\rho} d^\times t$ converges absolutely if $-(y - 2\rho)$ is an interior point of $\mathfrak{t}_{pc}^*$.

Let $t = (t_1, t_2)$ where

$$\begin{cases} t_1 = a_4(\underline{\lambda}_{11}, \underline{\lambda}_{12}, \underline{\lambda}_{13}, \underline{\lambda}_{14}), \ t_2 = a_2(\underline{\lambda}_{21}, \underline{\lambda}_{21}^{-1}) & \text{case (3)}, \\ t_1 = a_6(\underline{\lambda}_{11}, \underline{\lambda}_{12}, \underline{\lambda}_{13}, \underline{\lambda}_{14}, \underline{\lambda}_{15}, \underline{\lambda}_{16}), \ t_2 = a_2(\underline{\lambda}_{21}, \underline{\lambda}_{21}^{-1}) & \text{case (4)}. \end{cases}$$

Let $\gamma_{i,jk}$ be the weight of the coordinate $x_{i,jk}$ for all $i, j, k$ ($j > k$). The element $\gamma_{i,jk}$ can be expressed in the form $\gamma_{i,jk} = w(d_{i,jk})$, where $d_{i,jk} \in \mathbb{R}^4$ or $\mathbb{R}^6$.

Let $\sigma = \mathrm{Re}(s)$. We will prove that the function $\lambda^\sigma \Theta_{L_0}(\Psi, \underline{\lambda}t) t^{-2\rho}$ is integrable on $\mathbb{R}_+ \times T_\epsilon^0$ for $\sigma \gg 0$. What we are going to do is to divide $L_0$ into a union of finite number of (not necessarily $G_k$-stable) subsets $L_i$ and to estimate $\Theta_{L_i}(\Psi, \underline{\lambda}t) t^{-2\rho}$ by a finite number of functions of the form $\lambda^{p_N} t^{w(c_N)}$ where $p_N \in \mathbb{R}$, $c_N \in \mathbb{R}^4$ or $\mathbb{R}^6$ depend on a finite number of positive numbers $N$. These numbers should have the property that if we choose $N$ appropriately, $p_N \ll 0$ and all the entries of $c_N$ are negative.

If $\lambda \geq 1$, for any $\sigma \in \mathbb{R}$, we can choose $N$ depending on $\sigma$ so that $\sigma + p_N < 0$ and all the entries of $c_N$ are negative. This implies that the function $\lambda^\sigma \Theta_{L_i}(\Psi, \underline{\lambda}t) t^{-2\rho}$ is integrable on $[1, \infty) \times T_\epsilon^0$. If $\lambda \leq 1$, we fix $N$ so that all the entries of $c_N$ are negative. Then if $\sigma + p_N > 0$, the function $\lambda^\sigma \Theta_{L_i}(\Psi, \underline{\lambda}t) t^{-2\rho}$ is integrable on $(0, 1] \times T_\epsilon^0$. Since $\sigma$ is arbitrary for the convergence of the integral on $[1, \infty) \times T_\epsilon^0$, this proves the convergence of $Z(\Phi, s)$ for $\mathrm{Re}(s) \gg 0$ and $Z_+(\Phi, s)$ for all $s$.

Let

$$(3.3) \qquad I_0 = \begin{cases} \{(i, j, k) \mid i = 1, 2, \ 1 \leq k < j \leq 4\} & \text{case (3)}, \\ \{(i, j, k) \mid i = 1, 2, \ 1 \leq k < j \leq 6\} & \text{case (4)}. \end{cases}$$

For $I \subset I_0$, we define

$$(3.4) \qquad h_I(\lambda, t) = \prod_{(i,j,k) \in I} \sup(1, \lambda^{-1} t^{-\gamma_{i,jk}})$$

for $\lambda \in \mathbb{R}_+, t \in T_\epsilon^0$.

Functions of the form $h_I(\lambda, t)$ often appear in estimates of various incomplete theta series because our main tool is Lemma (1.2.6) [8]. So we first consider the function $h_I(\lambda, t)$. We start with the following two observations whose proofs are easy and are left to the reader.

**Lemma (3.5)**
(1) If $I_1 \subset I_2 \subset I_0$, then $h_{I_1}(\lambda, t) \leq h_{I_2}(\lambda, t)$.
(2) If $I = I_1 \coprod I_2 \subset I_0$, then $h_I(\lambda, t) = h_{I_1}(\lambda, t) h_{I_2}(\lambda, t)$.

**Lemma (3.6)**
$$h_I(\lambda, t) = \sup_{I' \subset I} \prod_{(i,j,k) \in I'} (\lambda^{-1} t^{-\gamma_{i,jk}}).$$

Next, to simplify the situation, we estimate $h_I(\lambda, t)$ by functions of the form $\lambda^p t^{w(c)}$.

Let

$$d_{i,jk} = \begin{cases} ((d_{i,jk,1}, \cdots, d_{i,jk,3}), d_{i,jk,4}) & \text{case (3)}, \\ ((d_{i,jk,1}, \cdots, d_{i,jk,5}), d_{i,jk,6}) & \text{case (4)}. \end{cases}$$

We define
$$c_{I,l} = \sum_{\substack{(i,j,k) \in I \\ d_{i,jk,l} < 0}} d_{i,jk,l},$$

for all $l$ and put
$$c_I = \begin{cases} ((c_{I,1}, \cdots, c_{I,3}), c_{I,4}) & \text{case (3)}, \\ ((c_{I,1}, \cdots, c_{I,5}), c_{I,6}) & \text{case (4)}. \end{cases}$$

**Lemma (3.7)** $h_I(\lambda, t) \ll \sup(1, \lambda^{-\#I}) t^{-w(c_I)}$ on $\mathbb{R}_+ \times T_\epsilon^0$.

*Proof.* Note that
$$\prod_{(i,j,k) \in I'} (\lambda^{-1} t^{-\gamma_{i,jk}}) = \lambda^{-\#I'} \prod_{(i,j,k) \in I'} t^{-w(d_{i,jk})}$$

and $\lambda^{-\#I'} \leq \sup(1, \lambda^{-\#I})$.

Let $\overline{d}_{i,jk} \in \mathbb{R}^4$ or $\mathbb{R}^6$ be the element obtained by replacing positive entries of $d_{i,jk}$ by 0. Then
$$\prod_{(i,j,k) \in I'} t^{-w(d_{i,jk})} \ll \prod_{(i,j,k) \in I'} t^{-w(\overline{d}_{i,jk})} = t^{-\sum_{(i,j,k) \in I'} w(\overline{d}_{i,jk})}.$$

However, since all the entries of $\overline{d}_{i,jk}$ are non-positive for all $i, j, k$.
$$t^{-\sum_{(i,j,k) \in I'} w(\overline{d}_{i,jk})} \ll t^{-\sum_{(i,j,k) \in I} w(\overline{d}_{i,jk})} = t^{-w(c_I)}.$$

This proves the lemma.
Q.E.D.

For the rest of this section, $\lambda \in \mathbb{R}_+$, $t \in T_\epsilon^0$. So in inequalities like Lemma (3.7), we will not mention that it is uniform with respect to $\lambda \in \mathbb{R}_+$, $t \in T_\epsilon^0$.

We first consider the case (3). The following lemma is easy to verify and the proof is left to the reader.

**Lemma (3.8)** (1) $d_{1,21} = ((\frac{1}{2}, 1, \frac{1}{2}), \frac{1}{2})$.
(2) $d_{1,31} = ((\frac{1}{2}, 0, \frac{1}{2}), \frac{1}{2})$.
(3) $d_{1,41} = ((\frac{1}{2}, 0, -\frac{1}{2}), \frac{1}{2})$.
(4) $d_{1,32} = ((-\frac{1}{2}, 0, \frac{1}{2}), \frac{1}{2})$.
(5) $d_{1,42} = ((-\frac{1}{2}, 0, -\frac{1}{2}), \frac{1}{2})$.
(6) $d_{1,43} = ((-\frac{1}{2}, -1, -\frac{1}{2}), \frac{1}{2})$.

Note that $d_{2,jk}$ can be obtained by replacing the last $\frac{1}{2}$ in $d_{1,jk}$ by $-\frac{1}{2}$. Also if we put $d_0 = ((-3, -4, -3), -1)$, then $-2\rho = w(d_0)$. By Lemma (3.7), $h_{I_0}(\lambda, t) \ll \sup(1, \lambda^{-12}) t^{-w(c_0)}$, where

(3.9) $$c_0 = c_{I_0} = -((3, 2, 3), 3).$$

**Definition (3.10)** (1) $L_1 = \{x \in L_0 \mid x_{1,21} \neq 0\}$.

(2) $L_2 = \{x \in L_0 \mid x_{1,21} = 0, \ x_{1,31} \neq 0\}$.
(3) $L_3 = \{x \in L_0 \mid x_{1,21} = 0, \ x_{1,31} = 0\}$.

Apparently, $L_0 = L_1 \cup L_2 \cup L_3$. So we estimate $\Theta_{L_i}(\Psi, \underline{\lambda}t)t^{-2\rho}$ for $i = 1, 2, 3$.

(1) Consider $L_1$.

Let $I = I_0 \setminus \{(1,2,1)\}$. By Lemma (1.2.6) [8], for any $N \geq 1$,

$$\Theta_{L_1}(\Psi, \underline{\lambda}t)t^{-2\rho} \ll \lambda^{-N} t^{-N\gamma_{1,21}} h_I(\lambda, t) t^{-2\rho}.$$

By Lemma (3.7),

$$h_I(\lambda, t) \ll \sup(1, \lambda^{-11}) t^{w(((3,2,3),3))},$$
$$h_I(\lambda, t) t^{-2\rho} \ll \sup(1, \lambda^{-11}) t^{w(((0,-2,0),2))}.$$

Since all the entries of $d_{1,21}$ are positive, $\lambda^\sigma \Theta_{L_1}(\Psi, \underline{\lambda}t) t^{-2\rho}$ is integrable on $\mathbb{R}_+ \times T_\epsilon^0$ for $\sigma \gg 0$ and on $[1, \infty) \times T_\epsilon^0$ for all $\sigma$.

(2) Consider $L_2$.

Let $I = I_0 \setminus \{(1,2,1), (1,3,1)\}$. By Lemma (1.2.6) [8], for any $N \geq 1$,

$$\Theta_{L_2}(\Psi, \underline{\lambda}t)t^{-2\rho} \ll \lambda^{-N} t^{-N\gamma_{1,31}} h_I(\lambda, t) t^{-2\rho}.$$

By Lemma (3.7),

$$h_2(\lambda, t) \ll \sup(1, \lambda^{-10}) t^{w(((3,2,3),3))},$$
$$h_2(\lambda, t) t^{-2\rho} \ll \sup(1, \lambda^{-10}) t^{w(((0,-2,0),2))}.$$

So for any $N \geq 1$,

$$\Theta_{L_2}(\Psi, \underline{\lambda}t) t^{-2\rho} \ll \lambda^{-N} \sup(1, \lambda^{-10}) t^{w(((0,-2,0),2) - N((\frac{1}{2},0,\frac{1}{2}),\frac{1}{2}))}.$$

Since all the entries of $((0,-2,0), 2) - N((\frac{1}{2}, 0, \frac{1}{2}), \frac{1}{2})$ are negative if $N > 4$, $\lambda^\sigma \Theta_{L_1}(\Psi, \underline{\lambda}t) t^{-2\rho}$ is integrable on $\mathbb{R}_+ \times T_\epsilon^0$ for $\sigma \gg 0$ and on $[1, \infty) \times T_\epsilon^0$ for all $\sigma$.

(3) Consider $L_3$.

Suppose $x \in L_3$. Then since $x_1$ is non-singular, $x_{1,32}, x_{1,41} \neq 0$. We define $I = I_0 \setminus \{(1,2,1), (1,3,1), (1,3,2), (1,4,1)\}$. Then by Lemma (1.2.6) [8], for any $N \geq 1$,

$$\Theta_{L_3}(\Psi, \underline{\lambda}t) t^{-2\rho} \ll \lambda^{-2N} t^{-N(\gamma_{1,32} + \gamma_{1,41})} h_I(\lambda, t).$$

By Lemma (3.7),

$$h_I(\lambda, t) \ll \sup(1, \lambda^{-8}) t^{w((\frac{5}{2}, 2, \frac{5}{2}), 3))},$$
$$h_I(\lambda, t) t^{-2\rho} \ll \sup(1, \lambda^{-8}) t^{w((-\frac{1}{2}, -2, -\frac{1}{2}), 2))}.$$

So for any $N \geq 1$,

$$\Theta_{L_3}(\Psi, \underline{\lambda}t) t^{-2\rho} \ll \lambda^{-2N} \sup(1, \lambda^{-8}) t^{w(((-\frac{1}{2}, -2, -\frac{1}{2}), 2) - N(d_{1,32} + d_{1,41}))}.$$

Since $d_{1,32} + d_{1,41} = ((0,0,0),1)$, all the entries of

$$\left(\left(-\frac{1}{2}, -2, -\frac{1}{2}\right), 2\right) - N(d_{1,32} + d_{1,41})$$

are negative if $N > 2$. Therefore, $\lambda^\sigma \Theta_{L_1}(\Psi, \underline{\lambda} t) t^{-2\rho}$ is integrable on $\mathbb{R}_+ \times T_\epsilon^0$ for $\sigma \gg 0$ and on $[1, \infty) \times T_\epsilon^0$ for all $\sigma$.

This completes the proof of Theorem (3.1) for the case (3).

Next, we consider the case (4). The following lemma is easy to verify and the proof is left to the reader.

**Lemma (3.11)** (1) $d_{1,21} = ((\frac{2}{3}, \frac{4}{3}, 1, \frac{2}{3}, \frac{1}{3}), \frac{1}{2})$.
(2) $d_{1,31} = ((\frac{2}{3}, \frac{1}{3}, 1, \frac{2}{3}, \frac{1}{3}), \frac{1}{2})$.
(3) $d_{1,41} = ((\frac{2}{3}, \frac{1}{3}, 0, \frac{2}{3}, \frac{1}{3}), \frac{1}{2})$.
(4) $d_{1,51} = ((\frac{2}{3}, \frac{1}{3}, 0, -\frac{1}{3}, \frac{1}{3}), \frac{1}{2})$.
(5) $d_{1,61} = ((\frac{2}{3}, \frac{1}{3}, 0, -\frac{1}{3}, -\frac{2}{3}), \frac{1}{2})$.
(6) $d_{1,32} = ((-\frac{1}{3}, \frac{1}{3}, 1, \frac{2}{3}, \frac{1}{3}), \frac{1}{2})$.
(7) $d_{1,42} = ((-\frac{1}{3}, \frac{1}{3}, 0, \frac{2}{3}, \frac{1}{3}), \frac{1}{2})$.
(8) $d_{1,52} = ((-\frac{1}{3}, \frac{1}{3}, 0, -\frac{1}{3}, \frac{1}{3}), \frac{1}{2})$.
(9) $d_{1,62} = ((-\frac{1}{3}, \frac{1}{3}, 0, -\frac{1}{3}, -\frac{2}{3}), \frac{1}{2})$.
(10) $d_{1,43} = ((-\frac{1}{3}, -\frac{2}{3}, 0, \frac{2}{3}, \frac{1}{3}), \frac{1}{2})$.
(11) $d_{1,53} = ((-\frac{1}{3}, -\frac{2}{3}, 0, -\frac{1}{3}, \frac{1}{3}), \frac{1}{2})$.
(12) $d_{1,63} = ((-\frac{1}{3}, -\frac{2}{3}, 0, -\frac{1}{3}, -\frac{2}{3}), \frac{1}{2})$.
(13) $d_{1,54} = ((-\frac{1}{3}, -\frac{2}{3}, -1, -\frac{1}{3}, \frac{1}{3}), \frac{1}{2})$.
(14) $d_{1,64} = ((-\frac{1}{3}, -\frac{2}{3}, -1, -\frac{1}{3}, -\frac{2}{3}), \frac{1}{2})$.
(15) $d_{1,65} = ((-\frac{1}{3}, -\frac{2}{3}, -1, -\frac{4}{3}, -\frac{2}{3}), \frac{1}{2})$.

Note that $d_{2,jk}$ can be obtained by replacing the last $\frac{1}{2}$ in $d_{1,jk}$ by $-\frac{1}{2}$. Also if we put $d_0 = ((-5, -8, -9, -8, -5), -1)$, then $-2\rho = w(d_0)$. By Lemma (3.7), $h_{I_0}(\lambda, t) \ll \sup(1, \lambda^{-30}) t^{-w(c_0)}$, where

$$(3.12) \qquad c_0 = c_{I_0} = -\left(\left(\frac{20}{3}, 8, 6, 8, \frac{20}{3}\right), \frac{15}{2}\right).$$

**Definition (3.13)** (1) $L_1 = \{x \in V_k^{\text{ss}} \mid x_{1,21} \text{ or } x_{1,31} \text{ or } x_{1,41} \neq 0\}$.
(2) $L_2 = \{x \in V_k^{\text{ss}} \mid x_{1,21} = x_{1,31} = x_{1,41} = 0, \ x_{1,32} \text{ or } x_{1,42} \neq 0\}$.
(3) $L_3 = \{x \in V_k^{\text{ss}} \mid x_{1,21} = x_{1,31} = x_{1,41} = x_{1,32} = x_{1,42} = 0\}$.
(4) $L_4 = \{x \in L_3 \mid x_{1,43}, x_{1,51} \neq 0\}$.
(5) $L_5 = \{x \in L_3 \mid x_{1,51} = 0, \ x_{1,43}, x_{1,52} \neq 0\}$.
(6) $L_6 = \{x \in L_3 \mid x_{1,51} = x_{1,52} = 0, \ x_{1,43}, x_{1,61} \neq 0\}$.
(7) $L_7 = \{x \in L_3 \mid x_{1,51} = x_{1,52} = x_{1,61} = 0, \ x_{1,43}, x_{1,62} \neq 0\}$.
(8) $L_8 = \{x \in L_3 \mid x_{1,51} = x_{1,52} = x_{1,61} = x_{1,62} = 0, \ x_{1,43}, x_{2,21} \neq 0\}$.
(9) $L_9 = \{x \in L_3 \mid x_{1,43} = 0, \ x_{1,51} \neq 0\}$.
(10) $L_{10} = \{x \in L_3 \mid x_{1,43} = x_{1,51} = 0, \ x_{1,52}, x_{1,61} \neq 0\}$.
(11) $L_{11} = \{x \in L_3 \mid x_{1,43} = x_{1,51} = x_{1,61} = 0, \ x_{1,52} \neq 0\}$.
(12) $L_{12} = \{x \in L_3 \mid x_{1,43} = x_{1,51} = x_{1,52} = 0, \ x_{1,61}, x_{1,53} \neq 0\}$.

(13) $L_{13} = \{x \in L_3 \mid x_{1,43} = x_{1,51} = x_{1,52} = x_{1,53} = 0, \ x_{1,61}, x_{1,54} \neq 0\}$.
(14) $L_{14} = \{x \in L_3 \mid x_{1,43} = x_{1,51} = x_{1,52} = x_{1,61} = 0\}$.
(15) $L_{15} = \{x \in L_{14} \mid x_{1,62}, x_{1,53} \neq 0\}$.
(16) $L_{16} = \{x \in L_{14} \mid x_{1,53} = 0, \ x_{1,62}, x_{1,54} \neq 0\}$.
(17) $L_{17} = \{x \in L_{14} \mid x_{1,62} = 0, \ x_{1,53} \neq 0\}$.
(18) $L_{18} = \{x \in L_{14} \mid x_{1,62} = x_{1,53} = 0, \ x_{1,63}, x_{1,54} \neq 0\}$.

**Proposition (3.14)** (1) $V_k^{\mathrm{ss}} = \coprod_{\substack{1 \leq i \leq 18 \\ i \neq 3, 14}} L_i$.

(2) If $x \in V_k^{\mathrm{ss}}$, there exist $1 \leq i \leq 2$, $2 \leq j \leq 6$ such that $x_{i,j1} \neq 0$.
(3) If $x \in L_6$, there exist $2 \leq j \leq 5$, $1 \leq k \leq 2$ $(j > k)$ such that $x_{2,jk} \neq 0$.
(4) If $x \in L_7$, there exists $2 \leq j \leq 5$ such that $x_{2,j1} \neq 0$.
(5) If $x \in L_9$ or $L_{10}$, there exist $1 \leq k < j \leq 4$ such that $x_{2,jk} \neq 0$.
(6) If $x \in L_{11}$, $x_{2,21}$ or $x_{2,31}$ or $x_{2,41} \neq 0$.
(7) If $x \in L_{12}$, there exist $2 \leq j \leq 4$, $1 \leq k \leq 2$ $(j > k)$ such that $x_{2,jk} \neq 0$.
(8) If $x \in L_{13}$, $x_{2,21}$ or $x_{2,31}$ or $x_{2,32} \neq 0$.
(9) If $x \in L_{15}$, $x_{2,21}$ or $x_{2,31}$ or $x_{2,41} \neq 0$.
(10) If $x \in L_{16}$, $x_{2,21}$ or $x_{2,31} \neq 0$.
(11) If $x \in L_{17}$ or $L_{18}$, $x_{2,21} \neq 0$.

*Proof.* Note that if $1 \leq i, j \leq 18$, $i, j \neq 3, 14$, and $i \neq j$, then $L_i \cap L_j = \emptyset$.

It is easy to see that if $x \in V_k^{\mathrm{ss}}$ and $x \notin L_1, L_2$, then $x \in L_3$. Suppose that $x \in L_3$ and $x_{1,43} \neq 0$. Then if $x \notin L_4, \cdots, L_7$, $x_{1,ij} = 0$ for $i = 1, 2$, $j = 2, \cdots, 6$. Suppose $x_{2,21} = 0$. Then by considering the cofactor expansion with respect to the first two columns, $\det M_x(v)$ is a product of $v_2^2$ and a sum of determinants of matrices of the form

$$v_1 \begin{bmatrix} 0 & 0 & 0 & 0 \\ 0 & 0 & 0 & 0 \\ * & * & * & * \\ * & * & * & * \end{bmatrix} + v_2 \begin{bmatrix} * & * & * & * \\ * & * & * & * \\ * & * & * & * \\ * & * & * & * \end{bmatrix}.$$

The determinant of the above matrix clearly is divisible by $v_2^2$. So $F_x(v)$ is divisible by $v_2^2$, which contradicts to the assumption $x \in V_k^{\mathrm{ss}}$. This implies that $x_{2,21} \neq 0$ and $x \in L_8$.

Suppose that $x \in L_3$ and $x_{1,43} = 0$. If $x_{1,51}$ or $x_{1,52} \neq 0$, then $x \in \coprod_{i=9}^{11} L_i$. So we assume that $x \in L_3$ and $x_{1,43} = x_{1,51} = x_{1,52} = 0$. Then $x_{1,61} \neq 0$ or $x \in L_{14}$. If $x_{1,61}, x_{1,53} \neq 0$, $x \in L_{12}$. Suppose $x_{1,61} \neq 0$, $x_{1,53} = 0$. Then if $x_{1,54} = 0$, $x_{1,jk} = 0$ for $j, k = 1, \cdots, 5$. So by the cofactor expansion with respect to the last row and the last column, $\det M_x(v)$ is divisible by $v_2^4$, which is a contradiction. Therefore, $x_{1,54} \neq 0$, which implies that $x \in L_{13}$.

Suppose $x \in L_{14}$. If $x_{1,62}, x_{1,53} \neq 0$, then $x \in L_{15}$. Also if $x_{1,62} = 0$, $x_{1,53} \neq 0$, then $x \in L_{17}$. Suppose $x_{1,62} \neq 0$, $x_{1,53} = 0$. If $x_{1,54} = 0$, then $x_{1,jk} = 0$ for $j, k = 1, \cdots, 5$, which cannot happen. So $x_{1,54} \neq 0$, which implies that $x \in L_{16}$. Suppose $x_{1,62} = x_{1,53} = 0$. Then $x_{1,54} \neq 0$ for the same reason. If $x_{1,63} = 0$, the first three columns of $x_1$ are zero. So $\det M_x(v)$ is divisible by $v_2^3$. But since $M_x(v)$ is an alternating matrix, $\det M_x(v)$ is divisible by $v_2^4$, which is a contradiction. This proves (1).

The statements (2), (3), (5) are clear.

Consider the statement (4). Let $x \in L_7$. Then the first column of $x_1$ is zero. Suppose $x_{2,j1} = 0$ for $j = 2, \cdots, 5$. Then by the cofactor expansion with respect to the $(6,1), (1,6)$-entries, $\det M_x(v)$ is a product of $x_{1,61}^2 v_2^2$ and the determinant of an alternating matrix of the form

$$v_1 \begin{bmatrix} 0 & 0 & 0 & 0 \\ 0 & 0 & * & * \\ 0 & * & 0 & * \\ 0 & * & * & 0 \end{bmatrix} + v_2 \begin{bmatrix} 0 & * & * & * \\ * & 0 & * & * \\ * & * & 0 & * \\ * & * & * & 0 \end{bmatrix}.$$

The determinant of the above matrix is divisible by $v_2^2$. So $F_x(v)$ is divisible by $v_2^2$, which is a contradiction. Therefore, if $x \in L_7$, there exists $2 \le j \le 5$ such that $x_{2,j1} \neq 0$.

Consider the statement (6). Suppose $x \in L_{11}$. If the statement of (6) is false, there exists $x \in L_3$ such that

(3.15) $\qquad x_{1,43} = x_{1,51} = x_{1,61} = x_{2,21} = x_{2,31} = x_{2,41} = 0,\ x_{1,52} \neq 0.$

We show that (3.15) cannot happen. Suppose (3.15) is satisfied. Then $M_x(v)$ is of the following form

$$v_1 \begin{bmatrix} 0 & 0 & 0 & 0 & 0 & 0 \\ 0 & 0 & 0 & 0 & * & * \\ 0 & 0 & 0 & 0 & * & * \\ 0 & 0 & 0 & 0 & * & * \\ 0 & * & * & * & 0 & * \\ 0 & * & * & * & * & 0 \end{bmatrix} + v_2 \begin{bmatrix} 0 & 0 & 0 & 0 & * & * \\ 0 & 0 & * & * & * & * \\ 0 & * & 0 & * & * & * \\ 0 & * & * & 0 & * & * \\ * & * & * & * & 0 & * \\ * & * & * & * & * & 0 \end{bmatrix}.$$

If $x_{2,51} = x_{2,61} = 0$, $\det M_x(v)$ is identically zero, which is a contradiction. So we assume that $x_{2,51}$ or $x_{2,61} \neq 0$, Let $g_1 \in \mathrm{GL}(6)_k$ be an element of the form

$$\begin{bmatrix} 1 & & & & & \\ & 1 & & & & \\ & & 1 & & & \\ & & & 1 & & \\ & & & & A & \end{bmatrix},$$

where $A \in \mathrm{GL}(2)_k$. By applying an element of the form $g = (g_1, I_2) \in \mathrm{GL}(6)_k \times \mathrm{GL}(2)_k$, we may assume that $x_{2,21} = x_{2,31} = x_{2,41} = x_{2,51} = 0$, $x_{2,61} \neq 0$. Note that by the action of $g$, $\det M_x(v)$ changes by a non-zero constant and the form of $x_1$ does not change.

Therefore, $\det M_x(v)$ is a product of $x_{2,61}^2 v_2^2$ and the determinant of an alternating matrix of the form

$$v_1 \begin{bmatrix} 0 & 0 & 0 & * \\ 0 & 0 & 0 & * \\ 0 & 0 & 0 & * \\ * & * & * & 0 \end{bmatrix} + v_2 \begin{bmatrix} 0 & * & * & * \\ * & 0 & * & * \\ * & * & 0 & * \\ * & * & * & 0 \end{bmatrix}.$$

The determinant of the above matrix is divisible by $v_2^2$. This implies that $F_x(v)$ is divisible by $v_2^2$, which is a contradiction.

Consider the statement (7). Let $x \in L_{12}$. If $x_{2,jk} = 0$ for $j = 2, 3, 4$, $k = 1, 2$, $M_x(v)$ is of the form

$$v_1 \begin{bmatrix} 0 & A_2 \\ A_1 & * \end{bmatrix} + v_1 \begin{bmatrix} 0 & B_2 \\ B_1 & * \end{bmatrix},$$

where $A_1, B_1$ are $2 \times 2$ and $A_2, B_2$ are $4 \times 4$. Also the first row of $A_1$, the first and the second columns of $A_2$ are zero. Since $\det M_x(v) = \det(v_1 A_1 + v_2 B_1) \det(v_1 A_2 + v_2 B_2)$, $\det M_x(v)$ is divisible by $v_2^3$. Since $M_x(v)$ is an alternating matrix, $\det M_x(v)$ is divisible by $v_2^4$, which is a contradiction. Therefore, there exists $1 \le k < j \le 4$ such that $x_{2,jk} \ne 0$.

Consider the statement (8). Let $x \in L_{13}$. If $x_{2,21} = x_{2,31} = x_{2,32} = 0$, $M_x(v)$ is of the form

$$v_1 \begin{bmatrix} 0 & -{}^t A \\ A & * \end{bmatrix} + v_1 \begin{bmatrix} 0 & -{}^t B \\ B & * \end{bmatrix},$$

where $A, B$ are $3 \times 3$ and the first and the second rows of $A$ are zero. Since $\det M_x(v) = \det(v_1 A + v_2 B)^2$, $\det M_x(v)$ is divisible by $v_2^4$, which is a contradiction. Therefore, $x_{2,21}$ or $x_{2,31}$ or $x_{2,32} \ne 0$.

Consider the statement (9). Let $x \in L_{15}$. Suppose $x_{2,21} = x_{2,31} = x_{2,41} = 0$. Then by the cofactor expansion with respect to the first row and the first column, $\det M_x(v)$ is a product of $v_2^2$ and the determinant of a matrix of the form

$$v_1 \begin{bmatrix} 0 & 0 & 0 & * \\ 0 & 0 & 0 & * \\ 0 & 0 & 0 & * \\ * & * & * & * \end{bmatrix} + v_2 \begin{bmatrix} * & * & * & * \\ * & * & * & * \\ * & * & * & * \\ * & * & * & * \end{bmatrix}.$$

Therefore, $\det M_x(v)$ is divisible by $v_2^4$, which is a contradiction. So $x_{2,21}$ or $x_{2,31}$ or $x_{2,41} \ne 0$.

Consider the statement (10). Let $x \in L_{16}$. Suppose $x_{2,21} = x_{2,31} = 0$. Then by the cofactor expansion with respect to the first row and the first column, $\det M_x(v)$ is a product of $v_2^2$ and the determinant of a matrix of the form

$$v_1 \begin{bmatrix} 0 & 0 & * & * \\ 0 & 0 & * & * \\ 0 & 0 & * & * \\ * & * & * & * \end{bmatrix} + v_2 \begin{bmatrix} * & * & * & * \\ * & * & * & * \\ * & * & * & * \\ * & * & * & * \end{bmatrix}.$$

Therefore, $\det M_x(v)$ is divisible by $v_2^3$. Since $M_x(v)$ is an alternating matrix, $\det M_x(v)$ is divisible by $v_2^4$, which is a contradiction. So $x_{2,21}$ or $x_{2,31} \ne 0$.

Consider the statement (11). Let $x \in L_{17}$ or $L_{18}$. Then $x_{1,62} = 0$. . Suppose $x_{2,21} = 0$. Then by considering the cofactor expansion with respect to the first row and the first column, $\det M_x(v)$ is a product of $v_2^2$ and the determinant of a matrix of the form

$$v_1 \begin{bmatrix} 0 & 0 & * & * \\ 0 & 0 & * & * \\ 0 & * & * & * \\ 0 & * & * & * \end{bmatrix} + v_2 \begin{bmatrix} * & * & * & * \\ * & * & * & * \\ * & * & * & * \\ * & * & * & * \end{bmatrix}.$$

So, $\det M_x(v)$ is divisible by $v_2^3$, which is a contradiction. Therefore, $x_{2,21} \neq 0$. This completes the proof of Proposition (3.14).

Q.E.D.

The following proposition is an immediate consequence of Proposition (3.14).

**Proposition (3.16)**

$$\Theta_{V_k^{ss}}(\Psi, \underline{\lambda}t) = \sum_{\substack{1 \leq i \leq 18 \\ i \neq 3, 14}} \Theta_{L_i}(\Psi, \underline{\lambda}t).$$

We now consider individual cases.

(1) Consider $L_1$.

Note that $t^{-\gamma_{i,21}}, t^{-\gamma_{i,31}} \ll t^{-\gamma_{i,41}}$. So by Lemma (1.2.6) [8] and Lemma (3.7), for any $N \geq 1$,

$$\Theta_{L_1}(\Psi, \underline{\lambda}t)t^{-2\rho} \ll \lambda^{-N} t^{-N\gamma_{1,41}} h_{I_0}(\lambda, t) t^{w(d_0)}$$
$$\ll \lambda^{-N} \sup(1, \lambda^{-30}) t^{w(d_0 - c_0 - Nd_{1,41})}.$$

Since

$$d_0 - c_0 - Nd_{1,41} = \left(\left(\frac{5}{3}, 0, -3, 0, \frac{5}{3}\right), \frac{15}{2}\right) - N\left(\left(\frac{2}{3}, \frac{1}{3}, 0, \frac{2}{3}, \frac{1}{3}\right), \frac{1}{2}\right),$$

all the entries of $-Nd_{1,41} - c_0 + d_0$ are negative if $N$ is large. Also if $N$ is large, the exponent of $\lambda$ tends to $-\infty$.

Therefore, $\lambda^\sigma \Theta_{L_1}(\Psi, \underline{\lambda}t)t^{-2\rho}$ is integrable on $\mathbb{R}_+ \times T_\epsilon^0$ for $\sigma \gg 0$ and on $[1, \infty) \times T_\epsilon^0$ for all $\sigma$.

(2) Consider $L_2$.

Let $I = I_0 \setminus \{(1,2,1), (1,3,1), (1,4,1)\}$. For $2 \leq \beta \leq 6$, we define

$$L_{2,1\beta} = \{x \in L_2 \mid x_{1,32}, x_{1,\beta 1} \neq 0\},$$
$$L_{2,2\beta} = \{x \in L_2 \mid x_{1,42}, x_{1,\beta 1} \neq 0\},$$
$$L_{2,3\beta} = \{x \in L_2 \mid x_{1,32}, x_{2,\beta 1} \neq 0\},$$
$$L_{2,4\beta} = \{x \in L_2 \mid x_{1,42}, x_{2,\beta 1} \neq 0\}.$$

(We only consider $\beta = 5, 6$ for $L_{2,1\beta}, L_{2,2\beta}$.)

By Proposition (3.14)(2), $L_2 = \cup_{\alpha, \beta} L_{2,\alpha\beta}$. So

$$\Theta_{L_2}(\Psi, \underline{\lambda}t) \leq \sum_{\alpha, \beta} \Theta_{L_{2,\alpha\beta}}(\Psi, \underline{\lambda}t).$$

We consider $L_{2,1\beta}$ first.

Let $I' = \{(1,3,2), (1,\beta,1)\}$, $I'' = I \setminus \{(1,3,2), (1,\beta,1)\}$. We define

$$V' = \{x \in V \mid x_{i,jk} = 0 \text{ for } (i,j,k) \notin I'\},$$
$$V'' = \{x \in V \mid x_{i,jk} = 0 \text{ for } (i,j,k) \notin I''\}.$$

The subsets $V', V''$ are subspaces of $V$, and $L_{2,1\beta}$ can be considered as a subset of $V'_k \oplus V''_k$. For $x \in V' \oplus V''$, let $p'(x), p''(x)$ be the projections to the first factor and the second factor respectively.

By Lemma (1.2.5) [8], there exist $0 \leq \Psi' \in \mathscr{S}(V'_{\mathbb{A}}), 0 \leq \Psi'' \in \mathscr{S}(V''_{\mathbb{A}})$ such that

$$\Theta_{L_{2,1\beta}}(\Psi, \underline{\lambda}t) \ll \sum_{x \in L_{2,1\beta}} \Psi'(p'(\underline{\lambda}tx))\Psi''(p''(\underline{\lambda}tx)).$$

We define

$$\Theta'(\Psi', \underline{\lambda}t) = \sum_{\substack{x \in V'_k \\ x_{1,32}, x_{1,\beta 1} \in k^\times}} \Psi'(\underline{\lambda}tx),$$

$$\Theta''(\Psi'', \underline{\lambda}t) = \sum_{x \in V''_k} \Psi''(\underline{\lambda}tx).$$

Then

$$\Theta_{L_{2,1\beta}}(\Psi, \underline{\lambda}t) \ll \Theta'(\Psi', \underline{\lambda}t)\Theta''(\Psi'', \underline{\lambda}t).$$

By Lemma (1.2.6) [8],

$$\Theta''(\Psi'', \underline{\lambda}t) \ll h_{I''}(\lambda, t) \leq h_I(\lambda, t).$$

We estimate $h_I(\lambda, t)$. We define

$$I_1 = \{(2, j, 1) \text{ for } j = 2, 3, 4\},$$
$$I_2 = \left\{ \begin{array}{l} (i, 5, 1), (i, 6, 1), (i, 3, 2), (i, 4, 2), (i, 5, 2), \\ (i, 6, 2), (i, 4, 3), (i, 5, 3), (i, 5, 4) \text{ for } i = 1, 2 \end{array} \right\},$$
$$I_3 = I \setminus (I_1 \cup I_2) = \{(i, 6, k) \text{ for } i = 1, 2, \ k = 3, 4, 5\}.$$

Then $I = I_1 \coprod I_2 \coprod I_3$. If $(i, j, k) \in I_1$, all the entries except for the last of $d_{i,jk}$ are positive. If $(i, j, k) \in I_2$, $d_{i,jk}$ is of the form $((c_1, *, *, *, -c_1), *)$ or $((*, c_2, *, -c_2, *), *)$.

By Lemma (3.7),

$$h_{I_1}(\lambda, t) \ll \sup(1, \lambda^{-3})t^{w(((0,0,0,0,0), \frac{3}{2}))},$$
$$h_{I_3}(\lambda, t) \ll \sup(1, \lambda^{-6})t^{w(((2,4,4,4,4), \frac{3}{2}))}.$$

We have to be a little more careful about $h_{I_2}(\lambda, t)$. By Lemma (3.6),

$$h_{I_2}(\lambda, t) \ll \sup(1, \lambda^{-18}) \sup_{I'_2 \subset I_2} \prod_{(i,j,k) \in I'_2} t^{-\gamma_{i,jk}}.$$

By the proof of Lemma (3.7), for each $I'_2 \subset I_2$, there exist $a, b \in \mathbb{R}$ such that

$$\prod_{(i,j,k) \in I'_2} t^{-\gamma_{i,jk}} \ll t^{w(((\frac{2}{3}, \frac{8}{3}, 2, \frac{4}{3}, \frac{4}{3}), \frac{9}{2}) + ((a,b,0,-b,-a),0))}.$$

So there exist a finite number of real numbers $a_1, \cdots, a_l, b_1, \cdots, b_l$ such that

$$h_{I_2}(\lambda, t) \ll \sup(1, \lambda^{-18}) \sum_{h=1}^{l} t^{w((\frac{2}{3}, \frac{8}{3}, 2, \frac{4}{3}, \frac{4}{3}), \frac{9}{2}) + ((a_h, b_h, 0, -b_h, -a_h))}.$$

Moreover, $-\frac{4}{3} \leq a_h \leq 4$, $-\frac{8}{3} \leq b_h \leq \frac{4}{3}$ for all $h$.

Let

(3.17) $$p_h = \left( \left( a_h - \frac{7}{3}, 0, -3, 0, -a_h + \frac{1}{3} \right), \frac{13}{2} \right),$$

$$q_h = \left( \left( a_h - \frac{7}{3}, b_h - \frac{4}{3}, -3, -b_h - \frac{8}{3}, -a_h + \frac{1}{3} \right), \frac{13}{2} \right).$$

Then we get the following lemma by the above considerations.

**Lemma (3.18)**

$$h_I(\lambda, t) t^{-2\rho} \ll \sup(1, \lambda^{-27}) \sum_{h=1}^{l} t^{w(q_h)} \ll \sup(1, \lambda^{-27}) \sum_{h=1}^{l} t^{w(p_h)}.$$

Therefore,

$$\Theta_{L_{2,1\beta}}(\Psi, \underline{\lambda}t) t^{-2\rho} \ll \sup(1, \lambda^{-27}) \sum_{h=1}^{l} \Theta'(\Psi', \underline{\lambda}t) t^{w(p_h)}.$$

This implies that we only have to estimate functions of the form $\Theta'(\Psi', \underline{\lambda}t) t^{w(p_h)}$.

For $L_{2,2\beta}, L_{2,3\beta}, L_{2,4\beta}$, exactly the same argument works replacing $I'$ by $I' = \{(1,4,2), (1,\beta,1)\}, \{(1,3,2), (2,\beta,1)\}, \{(1,4,2), (2,\beta,1)\}$ respectively.

By Lemma (1.2.6) [8], for any $N_1, N_2 \geq 1$,

$$\Theta'(\Psi', \underline{\lambda}t) t^{w(p_h)} \ll \lambda^{-N_1-N_2} \sup(1, \lambda^{-27}) t^{w(p_h - N_1 d_{1,32} - N_2 d_{1,\beta 1})}.$$

For $L_{2,2\beta}$, etc., we get the same estimate replacing $d_{1,32}$ or $d_{1,\beta 1}$ by $d_{1,42}$ or $d_{2,\beta 1}$. However, since $t^{-\gamma_{1,32}} \ll t^{-\gamma_{1,42}}$ and $t^{-\gamma_{1,\beta 1}}, t^{-\gamma_{2,\beta 1}} \ll t^{-\gamma_{2,61}}$ for all $\beta$, we only have to consider functions of the form

$$\lambda^{-N_1-N_2} \sup(1, \lambda^{-27}) t^{w(p_h - N_1 d_{1,42} - N_2 d_{2,61})}.$$

The point here is that we can choose $N_1, N_2$ for each $h$ separately. If we had used Lemma (1.2.6) [8] directly to $\Theta_{L_{2,1\beta}}(\Psi, \underline{\lambda}t)$, we get an estimate by the function

$$\lambda^{-N_1-N_2} \sup(1, \lambda^{-27}) t^{-N_1 \gamma_{1,32} - N_2 \gamma_{1,\beta 1}} \sum_{h=1}^{l} t^{w(p_h)},$$

and the choice of $N_1, N_2$ must be the same for all $h$. This is the reason why we had to separate the two non-zero coordinates to start with.

It is easy to see that

$$2d_{1,42} + d_{2,61} = \left((0,1,0,1,0), \frac{1}{2}\right).$$

We choose $N_1, N_2$ of the form $N_1 = 2N_3 + N_4, N_2 = N_3 + N_5$, where $N_3 \geq 1$, $N_4, N_5 \geq 0$. Then

$$p_h - N_1 d_{1,42} - N_2 d_{2,61}$$
$$= \left(\left(a_h - \frac{7}{3}, -N_3, -3, -N_3, -a_h + \frac{1}{3}\right), \frac{13 - N_3}{2}\right) - N_4 d_{1,42} - N_5 d_{2,61}.$$

If $a_h \geq 0$, we choose $N_4 = 4 + 3a_h, N_5 = 0$. Then

$$p_h - N_1 d_{1,42} - N_2 d_{2,61}$$
$$= \left(\left(-1, -\frac{4+3a_h}{3} - N_3, -3, -\frac{8+6a_h}{3} - N_3, -1\right), -\frac{4+3a_h}{2} + \frac{13-N_3}{2}\right).$$

Since $a_h \geq 0$,
$$-\frac{4+3a_h}{3}, -\frac{8+6a_h}{3}, -\frac{4+3a_h}{2} \leq 0.$$

Therefore,
$$t^{w(p_h - N_1 d_{1,42} - N_2 d_{2,61})} \ll t^{w(((-1,-N_3,-3,-N_3,-1), \frac{13-N_3}{2}))},$$
$$\lambda^{-N_4-N_5} \ll \sup(\lambda^{-4}, \lambda^{-8}) = \lambda^{-4} \sup(1, \lambda^{-4}).$$

If $a_h \leq 0$, we choose $N_4 = 4, N_5 = -\frac{3a_h}{2}$. Then

$$p_h - N_1 d_{1,42} - N_2 d_{2,61}$$
$$= \left(\left(-1, -\frac{4}{3} + \frac{a_h}{2} - N_3, -3, -\frac{8}{3} - \frac{a_h}{2} - N_3, -1\right), -2 + \frac{3a_h}{4} + \frac{13-N_3}{2}\right).$$

Since $-\frac{4}{3} \leq a_h$ for all $h$,
$$-\frac{4}{3} + \frac{a_h}{2}, \; -\frac{8}{3} - \frac{a_h}{2}, \; -2 + \frac{3a_h}{4} \leq 0.$$

Therefore,
$$t^{w(p_h - N_1 d_{1,42} - N_2 d_{2,61})} \ll t^{w(((-1,-N_3,-3,-N_3,-1), \frac{13-N_3}{2}))},$$
$$\lambda^{-N_4-N_5} \ll \sup(\lambda^{-4}, \lambda^{-10}) = \lambda^{-4} \sup(1, \lambda^{-6}).$$

By the above considerations,

$$\Theta'(\Psi', \underline{\lambda}t) t^{w(p_h)} \ll \sup(1, \lambda^{-33}) \lambda^{-3N_3-4} t^{w((-1,-N_3,-3,-N_3,-1), \frac{13-N_3}{2})}.$$

This bound does not depend on $\alpha, \beta, h$. So for any $N_3 \geq 1$,

$$\Theta_{L_2}(\Psi, \underline{\lambda}t)t^{-2\rho} \ll \sup(1, \lambda^{-33})\lambda^{-3N_3-4}t^{w((-1,-N_3,-3,-N_3,-1),\frac{13-N_3}{2})}.$$

Therefore, $\lambda^\sigma \Theta_{L_2}(\Psi, \underline{\lambda}t)t^{-2\rho}$ is integrable on $\mathbb{R}_+ \times T_\epsilon^0$ for $\sigma \gg 0$ and on $[1, \infty) \times T_\epsilon^0$ for all $\sigma$.

(3) Consider $L_4$.

Let $I = I_0 \setminus \{(1,2,1), \cdots, (1,5,1), (1,3,2), (1,4,2), (1,4,3)\}$. Then by Lemma (1.2.6) [8], for any $N \geq 1$,

$$\Theta_{L_4}(\Psi, \underline{\lambda}t)t^{-2\rho} \ll \lambda^{-3N}t^{-N(\gamma_{1,43}+2\gamma_{1,51})}h_I(\lambda, t)t^{-2\rho}.$$

By Lemma (3.7),

$$h_I(\lambda, t) \ll \sup(1, \lambda^{-23})t^{w((\frac{17}{3}, 8-\frac{2}{3}, 3, 8-\frac{1}{3}, \frac{20}{3}), \frac{15}{2})},$$
$$h_I(\lambda, t)t^{-2\rho} \ll \sup(1, \lambda^{-23})t^{w((\frac{2}{3}, -\frac{2}{3}, -3, -\frac{1}{3}, \frac{5}{3}), \frac{13}{2})}.$$

It is easy to see that

$$d_{1,43} + 2d_{1,51} = \left((1,0,0,0,1), \frac{3}{2}\right).$$

Since the first, the fifth, and the last entries are positive, all the entries of

$$\left(\left(\frac{2}{3}, -\frac{2}{3}, -3, -\frac{1}{3}, \frac{5}{3}\right), \frac{13}{2}\right) - N(d_{1,43} + 2d_{1,51})$$

are negative.

Therefore, $\lambda^\sigma \Theta_{L_4}(\Psi, \underline{\lambda}t)t^{-2\rho}$ is integrable on $\mathbb{R}_+ \times T_\epsilon^0$ for $\sigma \gg 0$ and on $[1, \infty) \times T_\epsilon^0$ for all $\sigma$.

(4) Consider $L_5$.

Let $L_{5,\alpha\beta} = \{x \in L_5 \mid x_{\alpha,\beta 1} \neq 0 \text{ for } \alpha = 1, 2, \beta = 2, \cdots, 6\}$. Then $L_5 = \cup_{\alpha,\beta} L_{5,\alpha\beta}$.

We define

$$I = I_0 \setminus \{(1,2,1), (1,3,1), (1,4,1), (1,5,1), (1,3,2), (1,4,2), (1,5,2)\}.$$

Then $h_I(\lambda, t)$ has the same bound as in Lemma (3.18) with $a_h, b_h, q_h \in \mathbb{R}$ for $h = 1, \cdots, l$.

We fix $\alpha, \beta$. Let $I' = \{(1,4,3), (1,5,2), (\alpha, \beta, 1)\}$. For $\Psi' \in \mathscr{S}(V'_\mathbb{A})$, we define

$$\Theta'(\Psi', \underline{\lambda}t) = \sum_{\substack{x \in V'_k \\ x_{1,43}, x_{1,52}, x_{\alpha,\beta,1}}} \Psi'(\underline{\lambda}tx).$$

By a similar consideration as before, there exists $0 \leq \Psi' \in \mathscr{S}(V'_\mathbb{A})$ such that

$$\Theta_{L_{5,\alpha\beta}}(\Psi, \underline{\lambda}t)t^{-2\rho} \ll \sup(1, \lambda^{-27}) \sum_{l=1}^{h} \Theta'(\Psi', \underline{\lambda}t)t^{w(q_h)}.$$

We consider each term. By Lemma (1.2.6) [8], for any $N_1, N_2, N_3 \geq 1$,

$$\Theta'(\Psi', \underline{\lambda}t)t^{w(q_h)} \ll \lambda^{-N_1-N_2-N_3} t^{w(q_h - N_1 d_{1,43} - N_2 d_{1,52} - N_3 d_{\alpha,\beta_1})}.$$

Since $t^{-\gamma_{\alpha,\beta_1}} \ll t^{-\gamma_{2,61}}$ for all $\alpha, \beta$ as above, we only consider the case $(\alpha, \beta) = (2, 6)$. It is easy to see that

$$d_{1,43} + d_{1,52} + d_{2,61} = \left((0,0,0,0,0), \frac{1}{2}\right).$$

So we put $N_1 = N_4+N_5$, $N_2 = N_4+N_6$, $N_3 = N_4+N_7$, where $N_4 \geq 1$, $N_5, N_6, N_7 \geq 0$.

Let $W = \{(a, b, 0, -b, -a) \mid a, b \in \mathbb{R}\} \subset \mathbb{R}^5$. The following lemma and its corollary are easy to verify and the proofs are left to the reader.

**Lemma (3.19)** *The convex hull of*

$$\left\{\left(-\frac{1}{3}, -\frac{2}{3}, 0, \frac{2}{3}, \frac{1}{3}\right), \left(-\frac{1}{3}, \frac{1}{3}, 0, -\frac{1}{3}, \frac{1}{3}\right), \left(\frac{2}{3}, \frac{1}{3}, 0, -\frac{1}{3}, -\frac{2}{3}\right)\right\}$$

*contains a neighborhood of the origin of $W$.*

**Corollary (3.20)** *For any $a, b \in \mathbb{R}$, there exist $c_1, c_2, c_3 \geq 0$ such that*

$$c_1\left(-\frac{1}{3}, -\frac{2}{3}, 0, \frac{2}{3}, \frac{1}{3}\right) + c_2\left(-\frac{1}{3}, \frac{1}{3}, 0, -\frac{1}{3}, \frac{1}{3}\right) + c_3\left(\frac{2}{3}, \frac{1}{3}, 0, -\frac{1}{3}, -\frac{2}{3}\right)$$
$$= (a, b, 0, -b, -a).$$

By the above Corollary, we choose $N_5, N_6, N_7 \geq 0$ so that

$$N_5\left(-\frac{1}{3}, -\frac{2}{3}, 0, \frac{2}{3}, \frac{1}{3}\right) + N_6\left(-\frac{1}{3}, \frac{1}{3}, 0, -\frac{1}{3}, \frac{1}{3}\right) + N_7\left(\frac{2}{3}, \frac{1}{3}, 0, -\frac{1}{3}, -\frac{2}{3}\right)$$
$$= \left(a_h - \frac{4}{3}, b_h, 0, -b_h, -a_h + \frac{4}{3}\right).$$

Then

$$\Theta'(\Psi', \underline{\lambda}t)t^{w(q_h)} \ll \lambda^{-3N_4-N_5-N_6-N_7} t^{w(((-1, -\frac{4}{3}, -3, -\frac{8}{3}, -1), \frac{13-N_4}{2}))}.$$

Since there are finitely many possibilities for $h$, there exist $c_1, c_2 \geq 0$ such that

$$\Theta_{L_{5,\alpha\beta}}(\Psi, \underline{\lambda}t)t^{-2\rho} \ll \lambda^{-3N_4-c_1} \sup(1, \lambda^{-27-c_2}) t^{w(((-1, -\frac{4}{3}, -3, -\frac{8}{3}, -1), \frac{13-N_4}{2}))}.$$

The right hand side does not depend on $\alpha, \beta, h$.

Therefore, $\lambda^\sigma \Theta_{L_5}(\Psi, \underline{\lambda}t)t^{-2\rho}$ is integrable on $\mathbb{R}_+ \times T_\epsilon^0$ for $\sigma \gg 0$ and on $[1, \infty) \times T_\epsilon^0$ for all $\sigma$.

(5) Consider $L_6$.

For $\alpha = 2, \cdots, 5$, $\beta = 1, 2$, $\alpha > \beta$, we define $L_{6,\alpha\beta} = \{x \in L_6 \mid x_{2,\alpha\beta} \neq 0\}$. Then by Proposition (3.14)(3), $L_6 = \cup_{\alpha,\beta} L_{6,\alpha\beta}$. Since $t^{-\gamma_{2,\alpha\beta}} \ll t^{-\gamma_{2,52}}$ for all $\alpha, \beta$ as above, we only consider the case $(\alpha, \beta) = (5, 2)$. If $x \in L_{6,52}$, $x_{1,43}, x_{1,61}, x_{2,52} \neq 0$. So by the same argument as in (4), $\lambda^\sigma \Theta_{L_6}(\Psi, \underline{\lambda}t)t^{-2\rho}$ is integrable on $\mathbb{R}_+ \times T_\epsilon^0$ for $\sigma \gg 0$ and on $[1, \infty) \times T_\epsilon^0$ for all $\sigma$.

(6) Consider $L_7$.

For $\alpha = 2, \cdots, 5$, we define $L_{7,\alpha} = \{x \in L_7 \mid x_{2,\alpha 1} \neq 0\}$. Then by Proposition (3.14)(4), $L_7 = \cup_\alpha L_{7,\alpha}$. Since $t^{-\gamma_{2,\alpha 1}} \ll t^{-\gamma_{2,51}}$ for all $\alpha$ as above, we only consider the case $\alpha = 5$

Let
$$I = I_0 \setminus \{(1, j, k) \text{ for } j = 2, \cdots, 6, \ k = 1, 2, \ j > k, \ (1, 4, 3)\}.$$

Then by Lemma (3.7),
$$h_I(\lambda, t) \ll \sup(1, \lambda^{-20}) t^{w(((5, 8-\frac{2}{3}, 6, 8-\frac{4}{3}, \frac{16}{3}), \frac{15}{2}))},$$
$$h_I(\lambda, t) t^{-2\rho} \ll \sup(1, \lambda^{-20}) t^{w(((0, -\frac{2}{3}, -3, -\frac{4}{3}, \frac{1}{3}), \frac{13}{2}))}.$$

By Lemma (1.2.6) [8], for any $N_1, N_2, N_3 \geq 1$,
$$\Theta_{L_{7,5}}(\Psi, \underline{\lambda}t) \ll \sup(1, \lambda^{-20}) \lambda^{-N_1 - N_2 - N_3}$$
$$\times t^{w(((0, -\frac{2}{3}, -3, -\frac{4}{3}, \frac{1}{3}), \frac{13}{2}) - N_1 d_{1,43} - N_2 d_{1,62} - N_3 d_{2,51})}.$$

It is easy to see that
$$3d_{1,43} + 2d_{1,62} + 4d_{2,51} = \left((1, 0, 0, 0, 1), \frac{1}{2}\right).$$

So if we choose $N_1 = 3N_4$, $N_2 = 2N_4$, $N_3 = 4N_4$ and $N_4 \gg 0$, all the entries of
$$\left(\left(0, -\frac{2}{3}, -3, -\frac{4}{3}, \frac{1}{3}\right), \frac{13}{2}\right) - N_1 d_{1,43} - N_2 d_{1,62} - N_3 d_{2,51}$$
$$= \left(\left(-N_4, -\frac{2}{3}, -3, -\frac{4}{3}, \frac{1}{3} - N_4\right), \frac{13 - N_4}{2}\right)$$

are negative.

Therefore, $\lambda^\sigma \Theta_{L_7}(\Psi, \underline{\lambda}t) t^{-2\rho}$ is integrable on $\mathbb{R}_+ \times T_\epsilon^0$ for $\sigma \gg 0$ and on $[1, \infty) \times T_\epsilon^0$ for all $\sigma$.

(7) Consider $L_8$.

By Lemma (1.2.6) [8], for any $N \geq 1$,
$$\Theta_{L_8}(\Psi, \underline{\lambda}t) t^{-2\rho} \ll \lambda^{-5N} t^{w(d_0 - N(3d_{1,43} + 2d_{2,21}))} h_{I_0}(\lambda, t).$$

It is easy to see that
$$3d_{1,43} + 2d_{2,21} = \left(\left(\frac{1}{3}, \frac{2}{3}, 2, \frac{10}{3}, \frac{5}{3}\right), \frac{1}{2}\right).$$

Since all the entries of the above element are positive, $\lambda^\sigma \Theta_{L_8}(\Psi, \underline{\lambda}t)t^{-2\rho}$ is integrable on $\mathbb{R}_+ \times T_\epsilon^0$ for $\sigma \gg 0$ and on $[1, \infty) \times T_\epsilon^0$ for all $\sigma$.

(8) Consider $L_9$.

Let
$$I = I_0 \setminus \{(1, j, k) \text{ for } 1 \leq k < j \leq 4, \ (1, 5, 1)\}.$$

Then by Lemma (3.7),
$$h_I(\lambda, t) \ll \sup(1, \lambda^{-23}) t^{w(((\frac{17}{3}, 8-\frac{2}{3}, 6, 8-\frac{1}{3}, \frac{20}{3}), \frac{15}{2}))},$$
$$h_I(\lambda, t) t^{-2\rho} \ll \sup(1, \lambda^{-23}) t^{w(((\frac{4}{3}, -\frac{2}{3}, -3, -\frac{1}{3}, \frac{5}{3}), \frac{13}{2}))}.$$

For $1 \leq \beta < \alpha \leq 4$, we define $L_{9,\alpha\beta} = \{x \in L_9 \mid x_{2,\alpha\beta} \neq 0\}$. Then by Proposition (3.14)(5), $L_9 = \cup_{\alpha,\beta} L_{9,\alpha\beta}$. Since $t^{-\gamma_{2,\alpha\beta}} \ll t^{-\gamma_{2,43}}$ for all $\alpha, \beta$ as above, by Lemma (1.2.6) [8], for any $N \geq 1$,

$$\Theta_{L_{9,\alpha\beta}}(\Psi, \underline{\lambda}t) t^{-2\rho} \ll \sup(1, \lambda^{-21}) \lambda^{-3N} t^{w(((0, -\frac{2}{3}, -3, -1, 1), \frac{13}{2}) - N(2d_{1,51} - d_{2,43}))}.$$

It is easy to see that
$$2d_{1,51} + d_{2,43} = \left((1, 0, 0, 0, 1), \frac{1}{2}\right).$$

So all the entries of
$$\left(\left(\frac{4}{3}, -\frac{2}{3}, -3, -\frac{1}{3}, \frac{5}{3}\right), \frac{13}{2}\right) - N(2d_{1,51} + d_{2,43})$$
$$= \left(\left(\frac{4}{3} - N, -\frac{2}{3}, -3, -\frac{1}{3}, \frac{4}{3} - N\right), \frac{13-N}{2}\right)$$

are negative if $N \gg 0$.

Therefore, $\lambda^\sigma \Theta_{L_9}(\Psi, \underline{\lambda}t) t^{-2\rho}$ is integrable on $\mathbb{R}_+ \times T_\epsilon^0$ for $\sigma \gg 0$ and on $[1, \infty) \times T_\epsilon^0$ for all $\sigma$.

(9) Consider $L_{10}$.

For $1 \leq \beta < \alpha \leq 4$, we define $L_{10,\alpha\beta} = \{x \in L_{10} \mid x_{2,\alpha\beta} \neq 0\}$. Then by Proposition (3.14)(5), $L_{10} = \cup_{\alpha,\beta} L_{10,\alpha\beta}$. Since $t^{-\gamma_{2,\alpha\beta}} \ll t^{-\gamma_{2,43}}$ for all $\alpha, \beta$ as above, we only consider the case $(\alpha, \beta) = (4, 3)$.

Let $I' = \{(1, 5, 2), (1, 6, 2), (2, 4, 3)\}$, and
$$V' = \{x \in V \mid x_{i,jk} = 0 \text{ for } (i, j, k) \notin I'\}.$$

For $\Psi' \in \mathscr{S}(V'_\mathbb{A})$, we define
$$\Theta'(\Psi', \underline{\lambda}t) = \sum_{\substack{x \in V'_k \\ x_{1,52}, x_{1,61}, x_{2,\alpha\beta} \neq 0}} \Psi'(\underline{\lambda}tx).$$

Then as before, there exists $0 \leq \Psi' \in \mathscr{S}(V'_{\mathbb{A}})$ such that

$$\Theta_{L_{10,43}}(\Psi, \underline{\lambda}t)t^{-2\rho} \ll \sup(1, \lambda^{-27})\Theta'(\Psi', \underline{\lambda}t) \sum_{h=1}^{l} t^{w(q_h)}.$$

By Lemma (1.2.6) [8], for any $N_1, N_2, N_3 \geq 1$,

$$\Theta'(\Psi', \underline{\lambda}t)t^{w(q_h)} \ll \sup(1, \lambda^{-27})\lambda^{-N_1-N_2-N_3}t^{w(q_h - N_1 d_{1,61} - N_2 d_{1,52} - N_3 d_{2,43})}.$$

It is easy to see that

$$d_{1,61} + d_{1,52} + d_{2,43} = \left((0, 0, 0, 0, 0), \frac{1}{2}\right).$$

Therefore, by the same argument as in (4), $\lambda^{\sigma}\Theta_{L_{10}}(\Psi, \underline{\lambda}t)t^{-2\rho}$ is integrable on $\mathbb{R}_+ \times T_{\epsilon}^0$ for $\sigma \gg 0$ and on $[1, \infty) \times T_{\epsilon}^0$ for all $\sigma$.

(10) Consider $L_{11}$.

For $\alpha = 2, 3, 4$, we define $L_{11,\alpha} = \{x \in L_{11} \mid x_{2,\alpha 1} \neq 0\}$. Then by Proposition (3.14)(6), $L_{11} = \cup_{\alpha} L_{11,\alpha}$. Since $t^{-\gamma_{2,\alpha 1}} \ll t^{-\gamma_{2,41}}$ for $\alpha = 2, 3, 4$, we only consider the case $\alpha = 4$.

It is easy to see that

$$3d_{1,52} + 2d_{1,41} = \left(\left(\frac{1}{3}, \frac{5}{3}, 0, \frac{1}{3}, \frac{5}{3}\right), \frac{1}{2}\right).$$

So by Lemma (1.2.6) [8], for any $N \geq 1$,

$$\Theta_{L_{11,4}}(\Psi, \underline{\lambda}t)t^{-2\rho} \ll \lambda^{-5N}t^{w(d_0 - N(3d_{1,52} + 2d_{1,41}))}h_{I_0}(\underline{\lambda}t)$$
$$\ll \lambda^{-5N}\sup(1, \lambda^{-30})t^{w(((\frac{5-N}{3}, -\frac{5N}{3}, -3, -\frac{N}{3}, \frac{5-5N}{3}), \frac{13-N}{2}))}.$$

Therefore, $\lambda^{\sigma}\Theta_{L_{11}}(\Psi, \underline{\lambda}t)t^{-2\rho}$ is integrable on $\mathbb{R}_+ \times T_{\epsilon}^0$ for $\sigma \gg 0$ and on $[1, \infty) \times T_{\epsilon}^0$ for all $\sigma$.

(11) Consider $L_{12}$.

For $\alpha = 2, 3, 4$, $\beta = 1, 2$, $\alpha > \beta$, we define $L_{12,\alpha\beta} = \{x \in L_{12} \mid x_{2,\alpha\beta} \neq 0\}$. Then by Proposition (3.14)(7), $L_{12} = \cup_{\alpha,\beta} L_{12,\alpha\beta}$. Since $t^{-\gamma_{2,\alpha\beta}} \ll t^{-\gamma_{2,42}}$ for all $\alpha, \beta$ as above, we only consider the case $(\alpha, \beta) = (4, 2)$.

It is easy to see that

$$3d_{1,61} + 2d_{1,53} + 4d_{2,42} = \left((0, 1, 0, 1, 0), \frac{1}{2}\right),$$
$$4d_{1,61} + 2d_{1,53} + 4d_{2,42} = \left(\left(\frac{2}{3}, \frac{4}{3}, 0, \frac{2}{3}, -\frac{2}{3}\right), 1\right),$$
$$d_{1,61} + d_{1,53} + 2d_{2,42} = \left(\left(-\frac{1}{3}, \frac{1}{3}, 0, \frac{2}{3}, \frac{1}{3}\right), 0\right).$$

Also
$$t^{-w(((\frac{2}{3},\frac{4}{3},0,\frac{2}{3},-\frac{2}{3}),1))} \ll t^{-w(((\frac{2}{3},0,0,0,-\frac{2}{3}),0))},$$
$$t^{-w(((-\frac{1}{3},\frac{1}{3},0,\frac{2}{3},\frac{1}{3}),0))} \ll t^{-w(((-\frac{1}{3},0,0,0,\frac{1}{3}),0))}.$$

We define
$$m_1 = \left((0,1,0,1,0), \frac{1}{2}\right),$$
$$m_2 = \left(\left(\frac{2}{3},0,0,0,-\frac{2}{3}\right), 0\right),$$
$$m_3 = \left(\left(-\frac{1}{3},0,0,0,\frac{1}{3}\right), 0\right).$$

Then by the same argument as in (2), we only have to consider functions of the form
$$\lambda^{-10N_1-4N_2-9N_3} \sup(1, \lambda^{-27}) t^{w(p_h - N_1 m_1 - N_2 m_2 - N_3 m_3)},$$
where $N_1 \geq 1, N_2, N_3 \geq 0$.

It is easy to see that we can choose $N_1, N_2, N_3$ so that all the entries of $p_h - N_1 m_1 - N_2 m_2 - N_3 m_3$ are negative.

Therefore, $\lambda^\sigma \Theta_{L_{12}}(\Psi, \underline{\lambda} t) t^{-2\rho}$ is integrable on $\mathbb{R}_+ \times T_\epsilon^0$ for $\sigma \gg 0$ and on $[1,\infty) \times T_\epsilon^0$ for all $\sigma$.

(12) Consider $L_{13}$.

For $(\alpha, \beta) = (2,1), (3,1), (3,2)$, we define $L_{13,\alpha\beta} = \{x \in L_{13} \mid x_{2,\alpha\beta} \neq 0\}$. Then by Proposition (3.14)(8) $L_{13} = \cup_{\alpha,\beta} L_{13,\alpha\beta}$. Since $t^{-\gamma_{2\alpha\beta}} \ll t^{-\gamma_{2,32}}$ for all $\alpha, \beta$ as above, we only consider the case $(\alpha, \beta) = (3, 2)$.

It is easy to see that
$$3d_{1,61} + 2d_{1,54} + 4d_{2,32} = \left((0,1,2,1,0), \frac{1}{2}\right),$$
$$d_{1,61} + d_{2,32} = \left(\left(\frac{1}{3}, \frac{2}{3}, 1, \frac{1}{3}, -\frac{1}{3}\right), 0\right),$$
$$d_{1,61} + d_{1,54} + 2d_{2,32} = \left(\left(-\frac{1}{3}, \frac{1}{3}, 1, \frac{2}{3}, \frac{1}{3}\right), 0\right).$$

Also
$$t^{-w(((\frac{1}{3},\frac{2}{3},1,\frac{1}{3},-\frac{1}{3}),0))} \ll t^{-w(((\frac{1}{3},0,0,0,-\frac{1}{3}),0))},$$
$$t^{-w(((-\frac{1}{3},\frac{1}{3},1,\frac{2}{3},\frac{1}{3}),0))} \ll t^{-w(((-\frac{1}{3},0,0,0,\frac{1}{3}),0))}.$$

Therefore, by the argument of (2) and (11), $\lambda^\sigma \Theta_{L_{13}}(\Psi, \underline{\lambda} t) t^{-2\rho}$ is integrable on $\mathbb{R}_+ \times T_\epsilon^0$ for $\sigma \gg 0$ and on $[1,\infty) \times T_\epsilon^0$ for all $\sigma$.

(13) Consider $L_{15}$.

For $\alpha = 2, 3, 4$, we define $L_{15,\alpha} = \{x \in L_{15} \mid x_{2,\alpha 1} \neq 0\}$. Then by Proposition (3.14)(9), $L_{15} = \cup_\alpha L_{15,\alpha}$. Since $t^{-\gamma_{2,\alpha 1}} \ll t^{-\gamma_{2,41}}$ for $\alpha = 2, 3, 4$, we only consider the case $\alpha = 4$.

It is easy to see that
$$2d_{1,62} + 2d_{1,53} + 3d_{2,41} = \left(\left(\frac{2}{3}, \frac{1}{3}, 0, \frac{2}{3}, \frac{1}{3}\right), \frac{1}{2}\right).$$

Then by Lemma (1.2.6) [8], for any $N \geq 1$,

$$\Theta_{L_{15,4}}(\Psi, \underline{\lambda}t)t^{-2\rho} \ll \lambda^{-7N}t^{w(d_0 - N((\frac{2}{3},\frac{1}{3},0,\frac{2}{3},\frac{1}{3}),\frac{1}{2}))}h_{I_0}(\lambda, t)$$
$$\ll \lambda^{-7N}\sup(1, \lambda^{-30})t^{w(((\frac{5}{3},0,-3,0,\frac{5}{3}),\frac{13}{2}) - N((\frac{2}{3},\frac{1}{3},0,\frac{2}{3},\frac{1}{3}),\frac{1}{2}))}.$$

Therefore, by the argument of (11), $\lambda^\sigma \Theta_{L_{15}}(\Psi, \underline{\lambda}t)t^{-2\rho}$ is integrable on $\mathbb{R}_+ \times T_\epsilon^0$ for $\sigma \gg 0$ and on $[1, \infty) \times T_\epsilon^0$ for all $\sigma$.

(14) Consider $L_{16}$.

For $\alpha = 2, 3$, we define $L_{16,\alpha} = \{x \in L_{16} \mid x_{2,\alpha 1} \neq 0\}$. Then by Proposition (3.14)(10), $L_{16} = \cup_\alpha L_{16,\alpha}$. Since $t^{-\gamma_{2,\alpha 1}} \ll t^{-\gamma_{2,31}}$ for $\alpha = 2, 3$, we only consider the case $\alpha = 3$.

It is easy to see that

$$2d_{1,62} + 2d_{1,54} + 3d_{2,31} = \left(\left(\frac{2}{3}, \frac{1}{3}, 1, \frac{2}{3}, \frac{1}{3}\right), \frac{1}{2}\right).$$

Since all the entries of the above element are positive, $\lambda^\sigma \Theta_{L_{16}}(\Psi, \underline{\lambda}t)t^{-2\rho}$ is integrable on $\mathbb{R}_+ \times T_\epsilon^0$ for $\sigma \gg 0$ and on $[1, \infty) \times T_\epsilon^0$ for all $\sigma$.

(15) Consider $L_{17}$.

It is easy to see that

$$3d_{1,53} + 2d_{2,21} = \left(\left(\frac{1}{3}, \frac{2}{3}, 2, \frac{1}{3}, \frac{5}{3}\right), \frac{1}{2}\right).$$

Since all the entries of the above element are positive, $\lambda^\sigma \Theta_{L_{17}}(\Psi, \underline{\lambda}t)t^{-2\rho}$ is integrable on $\mathbb{R}_+ \times T_\epsilon^0$ for $\sigma \gg 0$ and on $[1, \infty) \times T_\epsilon^0$ for all $\sigma$.

(16) Consider $L_{18}$.

It is easy to see that

$$2d_{1,63} + 2d_{1,54} + 3d_{2,21} = \left(\left(\frac{2}{3}, \frac{4}{3}, 1, \frac{2}{3}, \frac{1}{3}\right), \frac{1}{2}\right).$$

Since all the entries of the above element are positive, $\lambda^\sigma \Theta_{L_{18}}(\Psi, \underline{\lambda}t)t^{-2\rho}$ is integrable on $\mathbb{R}_+ \times T_\epsilon^0$ for $\sigma \gg 0$ and on $[1, \infty) \times T_\epsilon^0$ for all $\sigma$.

This completes the proof of Theorem (3.1) for the case (4).

Q.E.D.

Akihiko Yukie
Mathematics Department
College of Arts and Sceiences
Oklahoma State University
Stillwater OK 74078 USA